# CUQIpy: II. Computational uncertainty quantification for PDE-based inverse problems in Python

Amal M A Alghamdi[1], Nicolai A B Riis[1],
Babak M Afkham[1], Felipe Uribe[2,4],
Silja L Christensen[1], Per Christian Hansen[1]
and Jakob S Jørgensen[1,3,*]

[1] Department of Applied Mathematics and Computer Science, Technical University of Denmark. Richard Petersens Plads, Building 324, 2800 Kongens Lyngby, Denmark
[2] School of Engineering Sciences, Lappeenranta-Lahti University of Technology (LUT), Yliopistonkatu 34, 53850 Lappeenranta, Finland
[3] Department of Mathematics, The University of Manchester, Oxford Road, Alan Turing Building, Manchester M13 9PL, United Kingdom

E-mail: jakj@dtu.dk

## Abstract

Inverse problems, particularly those governed by Partial Differential Equations (PDEs), are prevalent in various scientific and engineering applications, and uncertainty quantification (UQ) of solutions to these problems is essential for informed decision-making. This second part of a two-paper series builds upon the foundation set by the first part, which introduced CUQIpy, a Python software package for computational UQ in inverse problems using a Bayesian framework. In this paper, we extend CUQIpy's capabilities to solve PDE-based Bayesian inverse problems through a general framework that allows the integration of PDEs in CUQIpy, whether expressed natively or using third-party libraries such as FEniCS. CUQIpy offers concise syntax that closely matches mathematical expressions, streamlining the modeling process and enhancing the user experience. The versatility and applicability of CUQIpy to PDE-based Bayesian inverse problems are demonstrated on examples covering parabolic, elliptic and hyperbolic PDEs. This includes problems involving the heat and Poisson equations and application case studies in electrical impedance tomography and photo-acoustic tomography, showcasing the software's efficiency, consistency, and intuitive interface. This comprehensive approach to UQ in PDE-based inverse problems provides accessibility for non-experts and advanced features for experts.

---

[4] Part of the work by F.U. was done while employed at the Technical University of Denmark.
[*] Author to whom any correspondence should be addressed.

## 1. Introduction

Inverse problems arise in various scientific and engineering applications, where the goal is to infer un-observable features from indirect observations. These problems are often ill-posed, making the inferred solution sensitive to noise in observed data and inaccuracies in forward models [19, 25]. Characterizing and evaluating uncertainties due to this sensitivity is crucial when making decisions based on inferred results.

To address these challenges, the field of Uncertainty Quantification (UQ) for inverse problems is in a phase of rapid growth [6, 39]. In medical imaging, for instance, UQ analysis allows experts to evaluate the uncertainty in cancer detection, which can directly impact patient treatment decisions [40]. In flood control and disaster management applications, UQ is needed to assess the risk of floods in specific regions, informing planning and resource allocation [30].

One particularly important category of inverse problems involves those governed by Partial Differential Equations (PDEs). These problems are encountered in various applications such as medical imaging [29, 41], seismic imaging [8, 9, 38], subsurface characterization [3, 4, 21], and non-destructive testing [34]. PDE-based inverse problems involve inferring parameters in PDE models from observed data, which introduces unique challenges in UQ due to the complex nature of the governing equations.

The Bayesian framework is widely used for UQ in both PDE-based and non-PDE-based inverse problems, as it enables the systematic incorporation of prior information, forward models, and observed data by characterizing the so-called posterior distribution [3, 4, 9, 11, 31]. This framework provides a comprehensive and unified approach for addressing the unique challenges of UQ in inverse problems.

### 1.1. Computational UQ for PDE-based inverse problems with CUQIpy

In this two-part series, we propose a Python software package CUQIpy, for Computational Uncertainty Quantification in Inverse problems. The package aims to make UQ analysis accessible to non-experts while still providing advanced features for UQ experts. The first paper [35] introduces the core library and components, and presents various test cases.

This second paper focuses on using CUQIpy (version 1.0.0) to solve Bayesian inverse problems where the forward models are governed by PDEs. While numerous software tools exist for modeling and solving PDE systems, such as FEniCS [28], FiPy [18], PyClaw [26], scikit-fem [17], and Firedrake [33], only few tools are specifically designed for PDE-based Bayesian inverse problems. The FEniCS-based package hIPPYlib [42, 43] is an example of a package that excels in this task.

To make UQ for PDE-based inverse problem more accessible, we propose a general framework for integrating PDE modeling tools into CUQIpy by defining an application programming interface (API) allowing PDE modeling libraries to interact with CUQIpy, regardless of the underlying PDE discretization scheme and implementation. This is possible because a major concept behind the design of **CUQIpy** is that the core components remain independent from specific forward modeling tools. On the other hand, plugins provide a flexible way to interface with third-party libraries, and in this paper we present **CUQIpy-FEniCS** as an example of a PDE-based plugin.

We introduce modules and classes in **CUQIpy** that enable solving PDE-based Bayesian inverse problems, such as `sampler`, `distribution`, and the `cuqi.pde` module. In the latter, the `cuqi.pde.PDE` class provides an abstract interface for integrating PDE modeling implementations like **FEniCS** with **CUQIpy**, simplifying the construction of PDE-based Bayesian inverse problems. The modules `cuqi.pde.geometry` and `cuqipy_fenics.pde.geometry` play an essential role, allowing the software to use information about the spaces on which the parameters and data are defined.

We demonstrate the versatility and applicability of **CUQIpy** through a variety of PDE-based examples, highlighting the integration and capabilities of the software. One example solves a Bayesian inverse problem governed by a one-dimensional (1D) heat equation, which underscores the intuitiveness of **CUQIpy**'s interface and its correspondence to the mathematical problem description. We present an elaborate electric impedance tomography (EIT) case study using the **CUQIpy-FEniCS** plugin, illustrating integration with third-party PDE modeling libraries. Finally, we examine a photoacoustic tomography (PAT) case, which shows **CUQIpy**'s ability to handle black-box forward models, emphasizing its adaptability to a wide range of applications in PDE-based Bayesian inverse problems. These examples effectively represent different classes of PDEs: parabolic, elliptic, and hyperbolic PDEs, respectively.

The examples involve inferring up to 100 parameters (the PAT problem) and solving PDEs of up to $4 \times 830$ state-variable dimensions (the EIT problem). They showcase different types of parameterization of the unknowns, namely, step expansion and Karhunen–Loève (KL) expansion and level set parameterization. The data setup varies in these examples from having data everywhere on the domain, having data only on parts of the domain or only on the boundaries, with noise levels explored of up to 20% of the noiseless data magnitude. The types of unknowns explored are spatially-varying PDE coefficients, as in the Poison and the EIT examples; and initial condition profiles, as in the heat and the PAT examples. Utilizing the flexibility of **CUQIpy** modeling framework, we demonstrate combining multiple datasets, namely, data resulting from multiple injection patterns in the EIT example. For simplicity, we assume these models are exact and we leave treatment of forward model error [7, 10] for future investigation. We emphasize that the **CUQIpy** framework is general to explore other types of unknown parameterizations, priors (e.g. using Markov random fields), and noise models; And model different unknown quantities, boundary conditions for example. Variations of these features are explored in the non-PDE-based Bayesian inverse problems presented in Part I [35] of this two-part series.

We have sought to design a versatile PDE abstraction layer for modeling a variety of PDE-based problems within the general Bayesian inverse problems framework provided by **CUQIpy** with focus on modularity and an intuitive, user-friendly interface. The goal of this paper is to demonstrate its utility on small to moderate scale problems on which we have found **CUQIpy** to perform well. Support for specialized PDE problems of high complexity and large-scale computing needs is an important area of development for the **CUQIpy** framework. We believe that the plugin structure, as exemplified by the **CUQIpy-FEniCS** plugin presented within, will provide a route to handle large-scale problems. This will combine the efficiency of dedicated third-party libraries (such as fluid dynamics solvers or implementations of adjoint equations, etc) with the convenience of the high-level modeling framework of **CUQIpy**.

*1.2. A motivating example*

We give a brief introductory example of the use of **CUQIpy** to solve a PDE-based inverse problem with the Poisson equation modeled in **FEniCS** using the **CUQIpy-FEniCS** plugin. More details of the underlying computational machinery are provided in section 3.

The inverse problem we consider is to infer a two-dimensional electric conductivity field $\sigma(\boldsymbol{\xi})$ of a unit square medium that lies in the domain $\bar{\Gamma} = [0,1]^2$, from a noisy observation of the electric potential measurement everywhere in the domain; we denote by $u(\boldsymbol{\xi})$ the electric potential and by $y(\boldsymbol{\xi})$ the observation of the electric potential, which in this case coincides with the solution $u$ on the entire domain, but in general may be available on a subset of the domain or a derived quantity.

The electric potential spatial distribution is governed by the 2D Poisson equation and is driven by a source term $f(\boldsymbol{\xi})$ and prescribed boundary conditions. The Poisson equation can be used to model other physical systems. For example, $\sigma$ can represent the thermal conductivity of a medium and $u$ its temperature; alternatively, $\sigma$ can represent the permeability of a porous medium and $u$ the pore-fluid pressure. The 2D Poisson equation we consider can be written as

$$\nabla \cdot \left( e^{w(\boldsymbol{\xi})} \nabla u(\boldsymbol{\xi}) \right) = f(\boldsymbol{\xi}) \qquad \text{for} \qquad \boldsymbol{\xi} \in \Gamma = (0,1)^2 \tag{1}$$

written here in terms of the log-conductivity field, i.e. $w(\boldsymbol{\xi}) = \log \sigma(\boldsymbol{\xi})$ to ensure positivity of the inferred conductivity field. In this example, we assume zero boundary conditions on the left and right boundaries of the square domain and zero Neumann boundary conditions on the top and bottom boundaries; and a source term $f(\boldsymbol{\xi}) = 1$. The forward problem concerns determining the observation $y(\boldsymbol{\xi})$ from a given log-conductivity $w(\boldsymbol{\xi})$. The inverse problem becomes the problem of inferring the log-conductivity $w(\boldsymbol{\xi})$ from an observed realization of $y(\boldsymbol{\xi})$.

In **CUQIpy** we consider the discretized form of this problem,

$$\boldsymbol{y} = \boldsymbol{A}(\boldsymbol{x}), \tag{2}$$

where $\boldsymbol{A}$ is a nonlinear forward model, which corresponds to solving the discretized PDE to produce the observation $\boldsymbol{y}$ from a log-conductivity given in terms of a parameter $\boldsymbol{x}$. **CUQIpy** (and in this case **CUQIpy-FEniCS**) provides a collection of demonstration test problems including one from which the present forward model can be obtained as:

```
A = FEniCSPoisson2D(dim=(32,32), field_type="KL", ...).model
```

Here, for brevity we have only shown a couple of the inputs to configure the problem. The PDE (1) is discretized using the finite-element method (FEM) and implemented using **FEniCS** on a structured triangular mesh on the physical domain $\Gamma$. The PDE solution and log-conductivity are approximated on a first-order Lagrange polynomial space, see, e.g. [14]. In this example, we consider the log-conductivity parameterized in terms of a truncated KL expansion [13] that enforces smoothness, to remedy the inherent instability of inferring coefficients of the Poisson equation [15]. The vector $\boldsymbol{x} = [x_1, \ldots, x_{n_{\text{KL}}}]^\top$ is the vector of expansion coefficients, here truncated at $n_{\text{KL}} = 32$.

In **CUQIpy** we consider $\boldsymbol{x}$ and $\boldsymbol{y}$ vector-valued random variables representing the parameter to be inferred and the data, respectively. To specify a Bayesian inverse problem, we express statistical assumptions on variables and the relations between them. Here, we assume an i.i.d. standard normal distribution on the KL expansion coefficients $\boldsymbol{x}$ and additive i.i.d. Gaussian noise with known standard deviation $s_{\text{noise}}$ on the data:

$$x \sim \text{Gaussian}\left(\mathbf{0}, \mathbf{I}\right) \tag{3a}$$

$$y \sim \text{Gaussian}\left(A\left(x\right), s_{\text{noise}}^2 \mathbf{I}\right). \tag{3b}$$

We can specify this in CUQIpy as (with np representing NumPy [20]):

```
x = Gaussian(np.zeros(n_KL), 1, geometry=G_KL)
y = Gaussian(A(x), s_noise**2, geometry=G_FEM)
```

We note the close similarity between the mathematical expressions and the syntax. Additionally the distributions have been equipped with so-called geometry object G_KL and G_FEM, which capture the interpretation of $x$ as KL coefficients and $y$ as FEM expansion coefficients; this is elaborated in section 3.

We consider a true log-conductivity as a sample from the prior distribution on $x$, which we conveniently generate and plot (figure 1(a)) by

```
x_true = x.sample()
x_true.plot()
```

where we note this is displayed as the log-conductivity FEM function, made possible by x being equipped with the G_KL geometry. The exact data $y^{\text{exact}}$ arising from $x^{\text{true}}$ can be determined (and plotted) as A(x_true).plot() while a noisy data realization $y^{\text{obs}}$ can be obtained by sampling $y$ conditioned on $x^{\text{true}}$ (figures 1(b) and (c)):

```
y_obs = y(x=x_true).sample()
y_obs.plot()
```

Again, knowledge of the geometry object G_FEM, in this case, enables visualizing $y^{\text{exact}}$ and $y^{\text{obs}}$ as FEM functions. CUQIpy provides a framework for specifying and solving Bayesian inverse problems through posterior MCMC sampling. In the most high-level case we simply specify a Bayesian Problem from the random variables $y$ and $x$, provide the observed data $y^{\text{obs}}$ and run the UQ() method:

```
BP = BayesianProblem(y, x).set_data(y=y_obs)
BP.UQ()
```

Under the hood, CUQIpy applies Bayes' theorem to construct the posterior distribution, selects a suitable sampler based on the problem structure (in this case the NUTS sampler [22]), samples the posterior and produces posterior mean and UQ plots (figure 1).

The results show that the mean is visually a reasonable approximation of the true conductivity. The variance magnitude is very small and tends to zero as $\xi$ gets closer to the left and right boundaries on which the PDE boundary conditions $u = 0$ are prescribed. Additionally, the computed credibility intervals (CIs) enclose the exact KL expansion coefficients. Approximately, the first 10 KL expansion coefficients are inferred with high certainty, and the general trend is that the uncertainty increases as the expansion mode number $i$ increases.

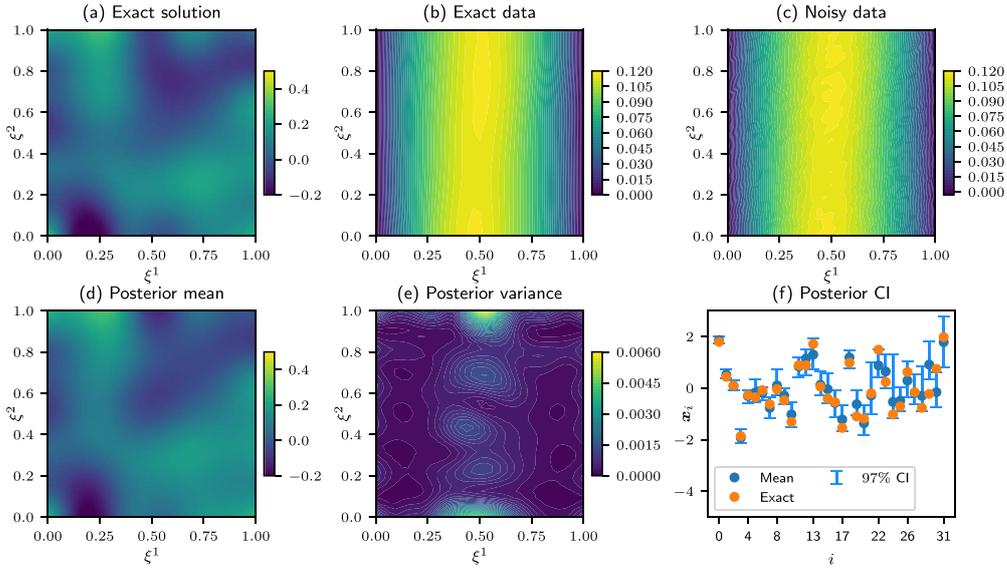

**Figure 1.** Results for the 2D Poisson problem, a prior sample $x^{\text{true}}$ is used as the exact solution. The noise level in the data is 1%. (a) The exact log-conductivity $w^{\text{true}} = G_{\text{KL}}(x^{\text{true}})$, see section 3 for $G_{\text{KL}}$. (b) The exact data $y^{\text{exact}}$. (c) The observed noisy data $y^{\text{obs}}$. (d) The log-conductivity mean: the posterior samples mean mapped through $G_{\text{KL}}$. (e) The log-conductivity variance computed from the posterior samples. (f) The CI plot showing the 97% CIs for the 32 KL coefficients (blue vertical lines), the exact KL coefficients $x^{\text{true}}$ (orange circles), and the KL coefficients means (blue circles).

### 1.3. Overview and notation

Having introduced and given a motivating example of UQ with CUQIpy for a PDE-based inverse problem in the present section, we present in section 2 our general framework for integrating PDE-based inverse problems in CUQIpy and illustrate this framework on an inverse problem governed by the 1D heat equation. In section 3, we describe our CUQIpy-FEniCS plugin that extends CUQIpy to allow UQ on PDE-based inverse problems modeled in FEniCS. We finish with two more elaborate case studies: First, in section 4 we demonstrate how electrical impedance tomography (EIT) with multiple layers of solution parametrization can be modeled with CUQIpy-FEniCS. Second, in section 5 we show how user-provided black-box PDE solvers can be used in CUQIpy in an example of Photo-Acoustic Tomography (PAT). Finally, in section 6 we conclude the paper.

We use the following notation: Calligraphic font such as $\mathcal{A}$ denotes a continuous operator. Bold upper case such as $\boldsymbol{A}$ denotes a discrete operator, with $\boldsymbol{I}_\ell$ denoting the $\ell \times \ell$ identity matrix; bold lower case such as $\boldsymbol{x}$ denotes a vector, and lower case such as $s$ and $f$ denotes a scalar or a scalar function with $p$ denoting a probability density function. We use the same notation for vectors and scalars to denote random vectors and scalars, to be distinguished by context. We denote by $\xi$ and $\boldsymbol{\xi} = [\xi^1, \xi^2]^{\mathsf{T}}$ the spatial coordinates in $\mathbb{R}$ and $\mathbb{R}^2$, respectively; and we denote by $\tau$ the time.

In the context of solving Bayesian inverse problems, we refer to the unknown quantity to be inferred as the *parameter* and the measured or observed quantities as the *data*, both considered random variables. When a superscript is provided for a parameter or a data vector, e.g. $x^{\text{true}}$, it indicates a particular realization of the parameter or the data, respectively. We refer to a particular noisy data realization that we use in the inversion, e.g. $y^{\text{obs}}$, as the *observed data*.

## 2. Framework for PDE-based Bayesian inverse problems in CUQIpy

In this section, we present our general framework for integrating PDE-based Bayesian inverse problems in CUQIpy. The framework is designed to be as generic as possible allowing—in principle—any PDE-based inverse problem to be handled. This includes PDEs expressed natively in CUQIpy (detailed in the present section), using a third-party PDE library such as FEniCS [28] (see sections 3 and 4) and through a user-provided black-box PDE-solver (see section 5). The critical components of this framework are provided by the `cuqi.pde` module, which contains the PDE class and its subclasses, the supporting Geometry classes, and the PDEModel class, see table 1.

The PDE class provides an abstract interface for representing PDEs, a subclass hereof is LinearPDE representing linear PDE problems. At present two concrete classes have been implemented: `SteadyStateLinearPDE` and `TimeDependentLinearPDE` from which a broad selection of steady-state and time-dependent linear PDEs can be handled. The Geometry classes allow us to parametrize in terms of various expansions to enforce desired properties on the solution, such as smoothness or as a step-wise function. The PDEModel class provides the interface to use the PDE as a CUQIpy Model for Bayesian inference. A PDEModel combines a PDE with two Geometry classes for the domain and range geometry to form a forward model of an inverse problem.

We illustrate this framework by an example: a Bayesian inverse problem governed by a 1D heat equation, sections 2.1–2.5. We emphasize that a much wider variety of PDEs can be handled; an example is used only for concreteness of the presentation.

### 2.1. The 1D heat equation inverse problem

We consider the inverse problem of reconstructing an initial temperature profile $g(\xi)$ at time $\tau = 0$ of a medium from temperature measurements $y(\xi)$ at a later time $\tau = \tau^{\max}$. We assume a medium that can be approximated by a 1D interval, $\xi \in [0,1]$. An example of such medium is a thin metal rod. The measurements are obtained over the interval $(0,1)$ or a subset of it, and they are typically polluted by a measurement error. The heat propagation in the medium from time $\tau = 0$ to time $\tau = \tau^{\max}$ can be modeled by a one-dimensional (1D) initial-boundary value heat equation, which can be written as:

$$\frac{\partial u(\xi,\tau)}{\partial \tau} - c^2 \frac{\partial^2 u(\xi,\tau)}{\partial \xi^2} = f(\xi,\tau), \quad \xi \in [0,1], \quad 0 \leqslant \tau \leqslant \tau^{\max}, \tag{4a}$$

$$u(0,\tau) = u(1,\tau) = 0, \tag{4b}$$

$$u(\xi,0) = g(\xi), \tag{4c}$$

where $u(\xi,\tau)$ is the temperature at time $\tau$ and location $\xi$, $c^2$ is the thermal conductivity (here taken to be a constant for simplicity), and $f$ is the source term. We assume zero boundary conditions, (4b), and an initial heat profile $g(\xi)$, (4c).

We define the *parameter-to-solution operator* $\mathcal{S}$ that maps the unknown parameter of the inverse problem $g(\xi)$ to the PDE solution $u(\xi,\tau)$, for $0 < \tau \leqslant \tau^{\max}$. Applying this operator is equivalent to solving the PDE (4a)–(4c) for a given initial condition $g(\xi)$. We also define the *observation operator* $\mathcal{O}$ that maps the PDE solution $u(\xi,\tau)$ to the observed quantities, the temperature measurements $y(\xi)$.

**Table 1.** A subset of CUQIpy classes that support integrating PDE-based problems. For a comprehensive list of classes and modules, see the companion paper [35].

| Class name | Description |
| --- | --- |
| `cuqi.pde` module: | |
| `PDE` | A class that represents the PDE, it implements the discretized maps $S$ and $O$ |
| `LinearPDE` | A class for linear PDE problems |
| `SteadyStateLinearPDE` | A class for steady-state linear PDE problems |
| `TimeDependentLinearPDE` | A class for time-dependent linear PDE problems |
| `cuqi.geometry` module: | |
| `Continuous1D` | A class that represents a 1D continuous space |
| `Continuous2D` | A class that represents a 2D continuous space |
| `KLExpansion` | A class for Karhunen–Loève expansion of functions |
| `StepExpansion` | A class for step functions |
| `cuqi.array` module: | |
| `CUQIarray` | A class for data arrays, subclassed from NumPy array |
| `cuqi.model` module: | |
| `PDEModel` | Forward model the uses a PDE-type class through calling its `assemble`, `solve` and `observe` methods |
| `cuqi.testproblem` module: | |
| `Poisson_1D` | 1D Poisson test problem (finite difference discretization) |
| `Heat_1D` | 1D Heat test problem (finite difference discretization) |

## 2.2. The discretized heat equation in CUQIpy

We discretize the system (4a)–(4c) in space using finite differences (FD). We discretize the solution $u(\xi,\tau)$ at a given time $\tau$ on a regular 1D grid of $n_{\text{grid}} = 100$ interior nodes. The grid spacing $h$ is approx. 0.01. We create a NumPy array to represent the grid

```
grid = np.linspace(h, 1-h, n_grid)
```

For simplicity, we use forward Euler for time stepping. For the choice $\tau^{\text{max}} = 0.01$, we discretize the time interval $[0,0.01]$ into $n_\tau = 225$ uniform time steps each of length $\Delta \tau$. We create a NumPy array to represent the time steps

```
tau = np.linspace(0, tau_max, n_tau)
```

We write the *k*th forward Euler step as follows

$$\boldsymbol{u}^{k+1} = \boldsymbol{u}^k + \Delta\tau\left(\boldsymbol{D}_c\boldsymbol{u}^k + \boldsymbol{f}^k\right), \quad \text{for } k = 0,\ldots,n_\tau, \tag{5}$$

where $\boldsymbol{u}^0 := \boldsymbol{g}$ is the initial condition $g$ discretized on the 1D grid, i.e. the $i$th element of $\boldsymbol{g}$ is $g(\xi_i)$ where $\xi_i$ is the coordinate of the $i$th grid node. Similarly, $\boldsymbol{u}^k$ and $\boldsymbol{f}^k$ are the PDE solution and the source term, respectively, at time $\tau = k\Delta\tau$ discretized on the 1D grid. $\boldsymbol{D}_c$ is the discretized diffusion operator $c^2 \partial^2/\partial\xi^2$, obtained using the centered-difference method. We create NumPy arrays to represent the right-hand side vector $\boldsymbol{f}^k$ (zero in this case) and the differential operator $\boldsymbol{D}_c$, and fix $c = 1$ for this example:

```
f = np.zeros(n_grid)
D_c = c**2 * ( np.diag(-2*np.ones(n_grid), 0) +
               np.diag(np.ones(n_grid-1), -1) +
               np.diag(np.ones(n_grid-1),  1) ) / h**2
```

We denote by $\boldsymbol{S}$ the discretized parameter-to-solution operator which maps the discretized initial condition $\boldsymbol{g}$ to the discretized PDE solution $\boldsymbol{u}$. $\boldsymbol{u}$ denotes the column vector of the time step solutions $\boldsymbol{u}^1, \ldots, \boldsymbol{u}^k, \ldots, \boldsymbol{u}^{n_\tau}$ stacked vertically. Additionally, we denote by $\boldsymbol{O}$ the discretized observation operator that maps the discretized PDE solution $\boldsymbol{u}$ to the observation $\boldsymbol{y} \in \mathbb{R}^m$, where $m$ is the number of measurements at locations $\{\xi_j^{\text{obs}}\}_{j=1}^m$. These locations might or might not correspond to the 1D grid points. In this example, they coincide with the set of grid points $\{\xi_i\}_{i=1}^{n_{\text{grid}}}$ or a subset of it.

To represent this discretized PDE equation in CUQIpy, we need to create a PDE-type object that encapsulates the details of the PDE equation and provides an implementation of the operators $\boldsymbol{S}$ and $\boldsymbol{O}$, table 1. Creating a PDE-type class requires a user-provided function that represents the components of the PDE on a standardized form, denoted by `PDE_form`. For time-dependent problems, the `PDE_form` function takes as inputs the unknown parameter that we want to infer (in this case g) and a scalar value for the current time: `tau_current`

```
def PDE_form(g, tau_current):
    return (D_c, f, g)
```

The `PDE_form` returns a tuple of the differential operator and right-hand side at time `tau_current` and the initial condition. Note that in this example, both the differential operator and the right-hand side, zero in this case, are independent of the time $\tau$.

For this 1D time-dependent heat equation, we create a `TimeDependentLinearPDE` object from the specific `PDE_form` and the time step vector `tau` and spatial grid `grid`:

```
PDE = TimeDependentLinearPDE(PDE_form, tau, grid_sol=grid)
```

The `TimeDependentLinearPDE` object calls the `PDE_form` every time step and passes the current time of the stepping method. The user can specify additional arguments when initializing the `TimeDependentLinearPDE` object, e.g. the spatial grid for observations, the time discretization scheme, and the linear solver to be used if the scheme is implicit. By default, the forward Euler method is used for time stepping and the observations are obtained at time $\tau^{\max}$ on the entire solution grid. We can print the PDE object using `print(PDE)`, which gives information about the object class and its `PDE_form`:

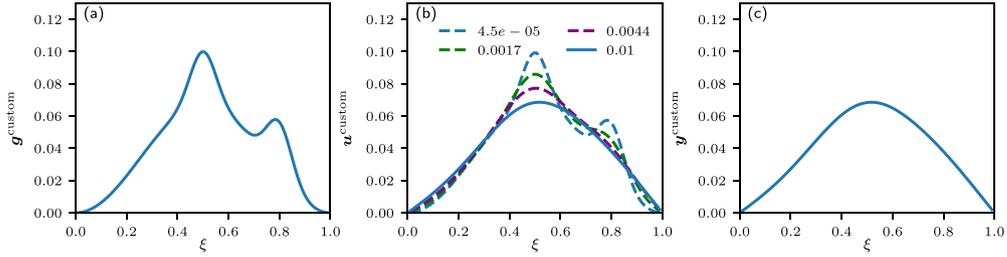

**Figure 2.** Illustration of the `TimeDependentLinearPDE` object's `assemble`, `solve` and `observe` methods. (a) Initial condition $g^{\text{custom}}$, (6), used as input to the `assemble` method. (b) PDE solution $u^{\text{custom}} = S(g^{\text{custom}})$, shown for selected times $\tau$ of the legend, obtained by the `solve` method. (c) Observation $y^{\text{custom}} = (O \circ S)(g^{\text{custom}})$, i.e. the PDE solution at time $\tau^{\max} = 0.01$ in this case, obtained by the `observe` method.

```
CUQI TimeDependentLinearPDE.
PDE form expression:
def PDE_form(g, tau_current):
    return (D_c, f, g)
```

All CUQIpy PDE-type classes implement three methods: (1) `assemble`, which performs any assembling that might be needed to prepare the matrices and vectors required to solve the PDE problem, (2) `solve`, which solves the PDE problem using the assembled components and is equivalent to applying the parameter-to-solution operator $S$, and (3) `observe`, which computes the observations from the PDE solution and is equivalent to applying the observation operator $O$. To illustrate these methods, let us consider an initial condition given by the expression

$$g^{\text{custom}}(\xi) = \frac{1}{30}\left(1 - \cos\left(2\pi\frac{1-\xi}{1}\right) + e^{-200(\xi-0.5)^2} + e^{-200(\xi-0.8)^2}\right). \quad (6)$$

We denote by $g^{\text{custom}}$ the discretization of $g^{\text{custom}}$ on the grid (see figure 2(a)). We call the method `assemble`, then apply the operator $S$ by calling the method `solve`:

```
PDE.assemble(g_custom)
u_custom, info = PDE.solve()
```

We show the solution `u_custom` in figure 2(b) where we plot selected time steps for illustration. Now we can apply the observation operator $O$, which in this case corresponds, conceptually, to a matrix that extracts the final time step solution $u^{n_\tau}$ from the entire PDE solution $u$. We denote the observation by $y^{\text{custom}} := u^{n_\tau}$ and show it in figure 2(c):

```
y_custom = PDE.observe(u_custom)
```

For time-dependent problems, PDE-type classes additionally implement the method `assemble_step` to assemble components that are needed to propagate the solution in time each time step, e.g. the discretized source term evaluated at time $\tau$. Furthermore, PDE-type classes can be equipped with the gradient of $O \circ S$ with respect to its input, $g$ in this case, in a given direction.

### 2.3. The 1D heat forward problem in CUQIpy

We define the discretized forward model of the 1D heat inverse problem

$$y = A(g) := (O \circ S)(g), \tag{7}$$

where $A : \mathbb{R}^n \to \mathbb{R}^m$. To represent the forward model in CUQIpy, we create an object from the class `PDEModel` which is a subclass of `Model`. To set up a `PDEModel` we need to specify which spaces are to be used for the domain and range of $A$; this is done using the `geometry` class. In the simplest case, the parameter $g$ and observation $y$ are simply vectors on the $\xi$ grid, which is specified by a `Continuous1D` geometry:

```
G_cont = Continuous1D(grid)
```

We can now set up the PDE model as

```
A = PDEModel(PDE, range_geometry=G_cont, domain_geometry=G_cont)
```

The `PDEModel` object encapsulates the PDE-type object and implements the `forward` method which corresponds to $A$. The `PDEModel` is agnostic to the underlying details of the PDE, e.g. the discretization method, the type of the PDE, and the third-party PDE modeling library used in the implementing the PDE Python methods. It uses the PDE object through calling the methods `assemble`, `solve`, and `observe`. One could continue with the present $A$ and solve directly for $g$, however here we demonstrate how to parametrize $g$ to enforce some desired properties on the inferred solution.

### 2.4. Parametrization by the `Geometry` class

The domain `Geometry` object represents the domain space of the forward model $A$. It can also be used to parametrize the unknown parameter, here $g$. As an example, we consider parameterization in terms of coefficients $x = [x_1, \ldots, x_{n_{\text{step}}}]^\top$ of an expansion

$$g = \sum_{i=1}^{n_{\text{step}}} x_i \chi_i, \tag{8}$$

where $\chi_i$ for $i = 1, \ldots, n_{\text{step}}$ is the characteristic function of the $i$th interval of a total of $n_{\text{step}}$ intervals in an equidistant partitioning of the domain $[0, 1]$. With this ('step expansion') parameterization of $g$, the unknown parameter of the inverse problem becomes the coefficients $x$. We denote by $G_{\text{step}}$ the discrete operator that maps $x$ to $g$. Thus we redefine the forward operator as

$$A(x) := (O \circ S \circ G_{\text{step}})(x), \tag{9}$$

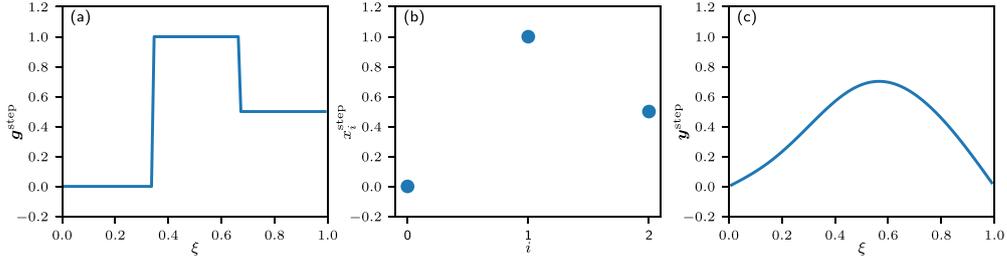

**Figure 3.** Illustration of CUQIpy `Geometry`. (a) `StepExpansion` function value plot of $g^{\text{step}} = G_{\text{step}}(x^{\text{step}})$. (b) `StepExpansion` parameter plot of $x^{\text{step}} = [0,1,0.5]^{\top}$. (c) `Continuous1D` function value plot of $y^{\text{step}} = A(x^{\text{step}})$, here with $\tau^{\max} = 0.02$.

where now $A : \mathbb{R}^{n_{\text{step}}} \to \mathbb{R}^m$. To specify the parameterization (8) in CUQIpy, we set up the domain geometry as a `StepExpansion` geometry object and pass the 1D grid and our choice of the number of steps `n_steps = 3` as arguments.

```
G_step = StepExpansion(grid, n_steps=3)
```

We can represent a function in this expansion by our fundamental data structure `CUQIarray`, which essentially contains a coefficient vector and the geometry, e.g.

```
x_step = CUQIarray([0, 1, 0.5], geometry=G_step)
```

`CUQIarray` has a dual representation: a *parameter value*, referring to the coefficient vector, here `[0, 1, 0.5]` and a *function value*, here the function with three steps considering parameters as expansion coefficients in the chosen geometry. A `Geometry`-type class implements the method `par2fun`, an implementation of the operator $G$, which maps the parameter value to the function value. It might also implement the method `fun2par`, the inverse map from the function to parameter value, $G^{-1}$, if it exists. It might also implement the gradient of $G$ with respect to $x$ in a given direction.

A `CUQIarray` allows convenient plotting of the object in context of the geometry:

```
x_step.plot()
```

By default, `plot` plots the function value representation of the variable, figure 3(a). That is, the call `x_step.plot()` results in calling the underlying `Geometry`-type object's `par2fun` method with the array values as the input and plotting its output, $g^{\text{step}} = G_{\text{step}}(x^{\text{step}})$. To plot the parameter value representation of the variable `x_step`, `plot_par = True` can be passed as an argument to the `plot` method, figure 3(b).

To employ the step function expansion we pass it as domain geometry:

```
A = PDEModel(PDE, range_geometry=G_cont, domain_geometry=G_step)
```

We can print the model A, `print(A)`, and get:

```
CUQI PDEModel: StepExpansion(3,) -> Continuous1D(100,).
    Forward parameters: ['x'].
    PDE: TimeDependentLinearPDE.
```

By default, the forward model input name is x. We can apply the forward model on $x^{\text{step}}$ and plot the result $y^{\text{step}} = A(x^{\text{step}})$:

```
y_step = A(x=x_step)
y_step.plot()
```

The returned `y_step` is a `CUQIarray` object equipped with the `G_cont` geometry (see figure 3(c)). Note, in this case we choose $\tau^{\max} = 0.02$, doubling the number of time steps.

*2.5. Specifying and solving the PDE-based Bayesian inverse problem*

In our discussion of the Bayesian modeling, we consider $x$ and $y$ to be random variables of the unknown parameter and the data, respectively. We are interested in a statistical characterization—the posterior distribution—of the unknown parameter $x$, given a prior distribution of $x$, a distribution of the data $y$, and a realization of the noisy data $y^{\text{obs}}$, see the companion paper for background on Bayesian modeling [35].

We define the Bayesian inverse problem, assuming additive Gaussian noise, as:

$$x \sim \text{Gaussian}\left(\mathbf{0}, \mathbf{I}_{n_{\text{step}}}\right),$$
$$y \sim \text{Gaussian}\left(A(x), s_{\text{noise}}^2 \mathbf{I}_m\right),$$

where $s_{\text{noise}}$ is the standard deviation of the data distribution, which we specify to dictate a desired noise level relative to the observed data, we assume a 10% noise level in this case. We use CUQIpy to create the distributions of $x$ and $y$ as follows

```
x = Gaussian(np.zeros(n_step), 1, geometry=G_step)
y = Gaussian(A(x), s_noise**2, geometry=G_cont)
```

We pass the argument `geometry = G_step` when initializing x to specify that samples from this distribution are expansion coefficients of the step expansion (8). Similarly, we pass the argument `geometry = G_cont` when initializing y. The argument `A(x)` represents that $y$ is conditioned on $x$ through the forward model, as shown by `print(y)`:

```
CUQI Gaussian. Conditioning variables ['x'].
```

We can draw five samples from the prior distribution and display these by:

```
prior_samples = x.sample(5)
prior_samples.plot()
```

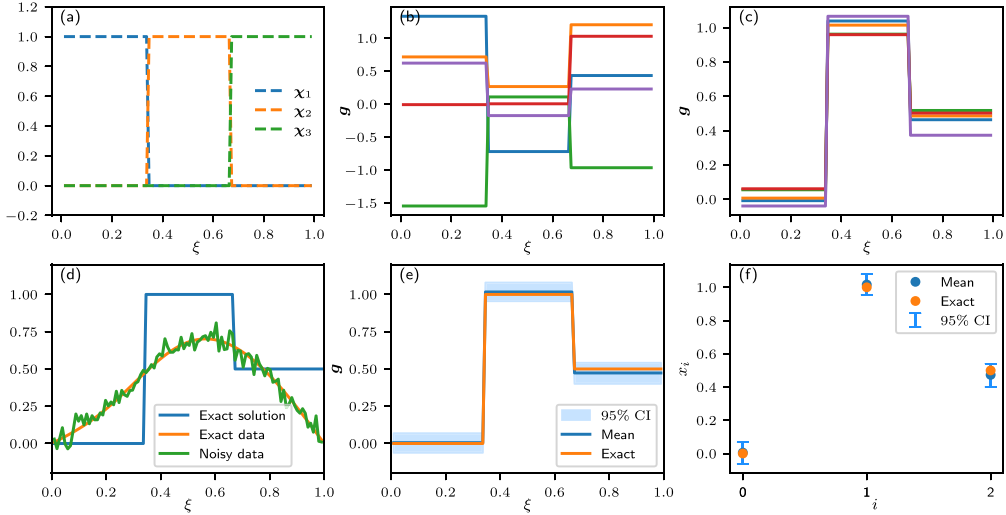

**Figure 4.** Results for the Bayesian inverse problem governed by the 1D heat equation in which we use the `StepExpansion` geometry and choose $\tau^{\max} = 0.02$. (a) Discretized characteristic functions $\boldsymbol{\chi}_1, \boldsymbol{\chi}_2,$ and $\boldsymbol{\chi}_3$, the basis functions of the expansion (8) for $n_{\text{step}} = 3$. (b) Prior samples plotted on the continuous domain. (c) Posterior samples plotted on the continuous domain. (d) The exact solution, exact data and noisy data. (e) The posterior sample mean and CI on the continuous domain. (f) The posterior sample means and CIs for the step expansion coefficients.

Here, `prior_samples` is a `Samples` object with the generated samples. It obtains the `Geometry`-type object from the prior x, here `G_step`. The prior samples are seen in figure 4(b). By default, the function values of the samples are plotted, i.e. the step functions.

We assume that the true solution is the step function with the coefficients $\boldsymbol{x}^{\text{step}}$, figure 3(b). We then generate synthetic noisy data $\boldsymbol{y}^{\text{obs}}$ by drawing a sample from the data distribution y conditioned on $\boldsymbol{x} = \boldsymbol{x}^{\text{step}}$:

```
y_obs = y(x=x_step).sample()
```

Figure 4(d) shows the exact solution $\boldsymbol{g}^{\text{step}}$, the exact data $\boldsymbol{y}^{\text{step}}$, and the noisy data $\boldsymbol{y}^{\text{obs}}$.

Now we have all the components we need to create the posterior distribution. We achieve this in CUQIpy by creating a joint distribution of the uncertain parameters $\boldsymbol{x}$ and the data $\boldsymbol{y}$ using the `JointDistribution` class, then we condition the joint distribution on the data $\boldsymbol{y}^{\text{obs}}$ to obtain the posterior distribution. The joint distribution is given by

$$p(\boldsymbol{x},\boldsymbol{y}) = p(\boldsymbol{y}|\boldsymbol{x})p(\boldsymbol{x}), \qquad (10)$$

where $p(\boldsymbol{x})$ is the prior probability density function (PDF) and $p(\boldsymbol{y}|\boldsymbol{x})$ is the data distribution PDF. In CUQIpy, this translates to

```
joint = JointDistribution(x, y)
posterior = joint(y=y_obs) # condition on y=y_obs
```

Calling `print(joint)`, for example, gives:

```
JointDistribution(
    Equation:
        p(x,y) = p(x)p(y|x)
    Densities:
        x ~ CUQI Gaussian.
        y ~ CUQI Gaussian. Conditioning variables ['x'].
)
```

**CUQIpy** uses MCMC sampling methods, provided by its `Sampler` classes, to approximate the posterior and compute its moments, in particular, mean and variance. In this example, we use a component-wise Metropolis Hastings (CWMH) sampler [35, § 2] and set up an instance of it by simply passing the posterior as input:

```
my_sampler = CWMH(posterior)
```

`Sampler`-type classes implement the methods `sample` and `sample_adapt`. The latter adjusts the sampling scale (step size) to achieve a target acceptance rate, which is method dependent. For the `CWMH` sampler, the target acceptance rate is approx. 23%.

We generate 50 000 samples using the `CWMH` sampler:

```
posterior_samples = my_sampler.sample_adapt(50000)
```

`posterior_samples` is a `Samples` object which contains, in addition to the samples and their corresponding geometry object, the sampling acceptance and rejection information.

### 2.6. Posterior samples analysis, and visualization

The `Samples` class provides analysis, and visualization methods that can be used to study the posterior samples. Some of these methods integrate functionalities from **ArviZ**, a python package for exploratory analysis of Bayesian models [27]. For brevity, we only show some of the visualization features and refer the reader to **CUQIpy**'s documentation for more information on visualization.

A basic `Samples` operation is to plot selected samples (figure 4(c)):

```
posterior_samples.plot([2000, 3000, 4000, 5000, 6000])
```

We visualize the samples credibility interval (CI) using the method `plot_ci` which generates a plot of the samples CI, the sample mean, and the exact solution of the Bayesian inverse problem, if provided:

```
posterior_samples.plot_ci(95, exact=x_step)
```

The first argument is the CI expressed in percent, 95% CI in this case, and the second optional argument is the exact solution. In figure 4(e), we show the CI plot. Note that in this plot, the CI is plotted over the continuous domain $(0, 1)$ and that the CI encloses the exact solution. We can alternatively plot the CI for the coefficients $x_i$ by passing the argument `plot_par = True` to the `plot_ci` function, see figure 4(f) for the coefficient CI plot. In this coefficient CI plots, we also note that $x^{\text{step}}$ lies within the CI.

*2.7. Parameterizing the initial condition using KL expansion*

Here we present a different parameterization of the unknown initial condition $g(\xi)$ to elaborate on CUQIpy's modeling capabilities; we use a truncated KL expansion [23, 44]. Using this representation we are able to impose some regularity and spatial correlation on $g(\xi)$ and reduce the dimension of the discretized unknown parameter from $n$ to $n_{\text{KL}}$, where $n_{\text{KL}} \ll n$.

To do this we wish to express our $g$ as a vector-valued random variable following a zero-mean Gaussian distribution with a carefully constructed covariance matrix $C$ capturing the desired variance and spatial correlation. In this particular case, $C$ will be constructed as $C = \frac{1}{a^2} E \Lambda^2 E^T$, where the matrix $E$ is an $n \times n_{\text{KL}}$ matrix with ortho-normal columns, $1/a^2$ is the variance, and $\Lambda$ is an $n_{\text{KL}} \times n_{\text{KL}}$ diagonal matrix with diagonal elements $\lambda_i = 1/i^\gamma$, where $\gamma$ is a constant that controls the decay rate of the diagonal elements. The columns of $E$ are often chosen to be a discretization of continuous functions on a grid. Here, we choose the sinusoidal basis functions. This choice ensures that the boundary condition (4*b*) is imposed on the initial condition $g$. It can be shown that $g$ follows the desired distribution if we express it as

$$g = \frac{1}{a} \sum_{i=1}^{n_{\text{KL}}} x_i \sqrt{\lambda_i} e_i. \tag{11}$$

Here, $x_i$, $i = 1, \ldots, n_{\text{KL}}$ are independent standard normal random variables, known as KL expansion coefficients, and $e_i$ is the *i*th column of $E$. We show a few basis functions $e_i$ discretized on `grid` in figure 5(a). The expansion in (11) is known as the KL-expansion and, if $n_{\text{KL}} < n$, the expansion is truncated and $n_{\text{KL}}$ is the truncation size. This parameterization is a suitable choice for inferring the initial condition in the heat equation because the corresponding forward model $S$ is a smoothing operator with rapid eigenvalues decay, namely, exponential decay [15]. Thus, without parameterization or regularization, recovering the oscillatory components of the initial condition is unstable.

We denote by $G_{\text{KL}}$ the operator which maps the KL expansion coefficients vector $x = [x_1, \ldots, x_{n_{\text{KL}}}]$ to the approximated discretized initial condition $g$. We set up the domain geometry as a `KLExpansion` geometry and pass the arguments `decay_rate = 1.5`, `normalizer = 10`, and `num_modes = 20` for $\gamma$, $a$ and $n_{\text{KL}}$, respectively:

```
G_KL = KLExpansion(grid, decay_rate=1.5, normalizer=10, num_modes=20)
```

As in case of the step expansion, we then set up the prior as a `Gaussian` distribution with zero mean and identity covariance, passing also the argument `geometry = G_KL` and sample the prior and plot its samples, figure 5(b).

We use the custom initial condition $g^{\text{custom}}$ in (6) as the true solution. Then following the steps in section 2.5, we create the corresponding synthesized data $y^{\text{obs}}$. We study three cases, using this initial condition: 0.1% noise case (figure 5, second row), 5% noise case (figure 5, third row), and 5% noise with data available only on the first half of the domain (figure 5, fourth row). In the first two cases, we have data measurement everywhere in the domain. To specify the limited observation in the third case, we pass `grid_obs = grid[:50]` to the `TimeDependentLinearPDE` initializer. We also pass `grid[:50]` when creating the range geometry, instead of passing the entire `grid`. We then create the posterior, and sample it using the `CWMH` sampler, posterior samples of the second case are shown in figure 5(c). We note that as the noise level increases, the width of the continuous CI increases and less modes are reconstructed with high certainty. Also, observing only on the first half of the domain leads to a significantly wider CI in the part of the domain where we do not have data, figure 5(k), and higher uncertainties in the mode reconstructions, figure 5(l).

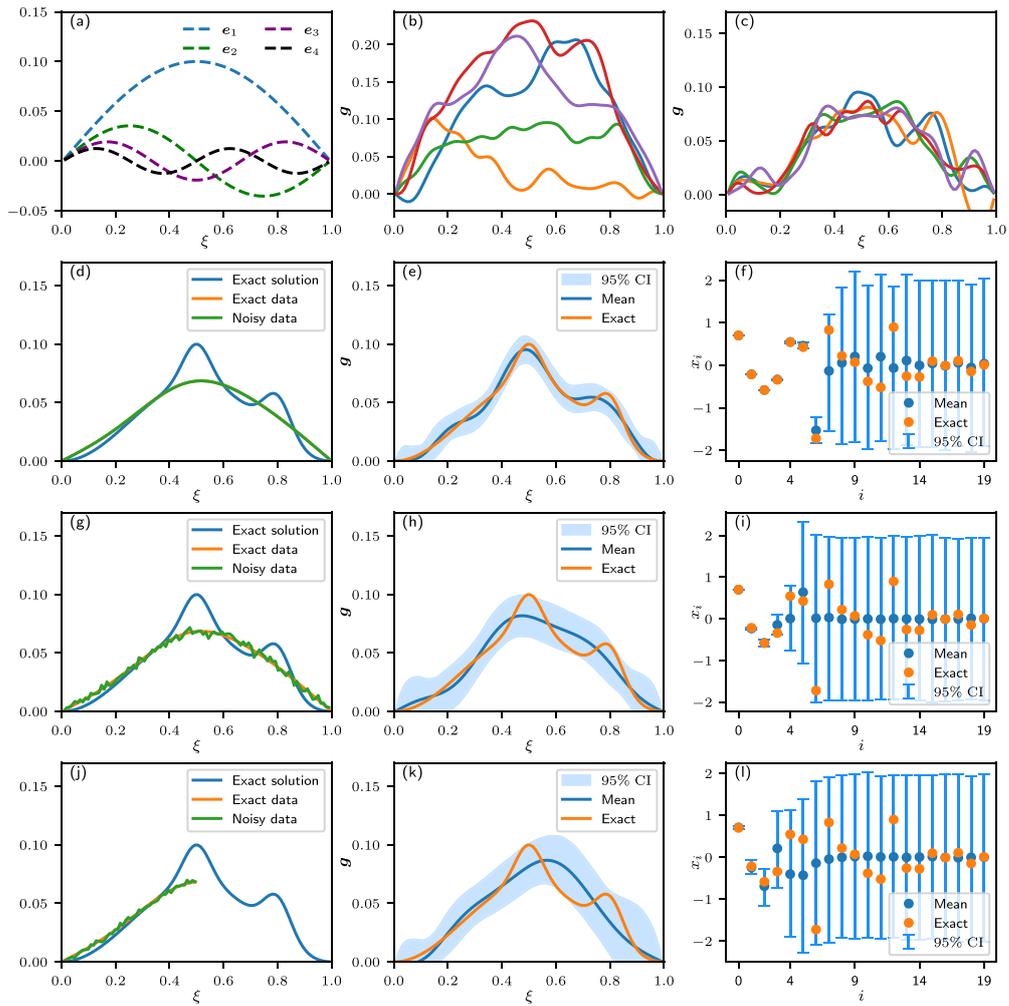

**Figure 5.** Results for the 1D heat equation-based Bayesian inverse problem in which we use the `KLExpansion` geometry and the function $g^{\text{custom}}$ (6) as the exact solution; and set $\tau^{\max} = 0.01$. We study three cases: 0.1% noise level case (second row), 5% noise level case (third row), and 5% noise level case and data available on the interval $(0, 0.5)$ only (fourth row). For the first two cases the data is available everywhere in the domain. (a) The KL expansion (11) basis functions $e_i$, for $i = 1, 2, 3, 4$. (b) Prior samples plotted on the continuous domain. (c) Posterior samples plotted on the continuous domain for the second case. For each case, the first column shows the exact solution $g^{\text{custom}}$, the exact data $y^{\text{custom}}$ and the observed noisy data $y^{\text{obs}}$, the second column shows the posterior sample mean and CI on the continuous domain, and the third column shows the posterior sample means and CIs for the KL expansion coefficients.

This concludes our overview of solving PDE-based Bayesian inverse problems with CUQIpy. We emphasize that the heat equation example was for demonstration and that the framework can be applied to wide variety of PDE-based inverse problems. In the next section, we show how to handle problems modeled in the FEM platform FEniCS.

## 3. CUQIpy-FEniCS plugin

FEniCS [28] is a popular Python package for solving PDEs using the FEM. The extent of its users, both in the academia and the industry, makes a dedicated CUQIpy interface for FEniCS highly desirable. Here we present our interface plugin: CUQIpy-FEniCS, and revisit the 2D Poisson example discussed in section 1.2 to unpack the underlying CUQIpy and CUQIpy-FEniCS components that are used in building the example. In section 4, we present an elaborate test case of using CUQIpy together with CUQIpy-FEniCS to solve an EIT problem with multiple data sets and in section 5 we showcase using some of the CUQIpy-FEniCS features in solving a PAT problem with a user-provided forward model.

We use the main modules of FEniCS: `ufl`, the FEniCS unified form language module, and `dolfin`, the Python interface of the computational high-performance FEniCS C++ backend, DOLFIN. We can import these modules as follows:

```
import ufl
import dolfin as dl
```

The CUQIpy-FEniCS plugin structure can be adopted to create CUQIpy plugins integrating other PDE-based modeling libraries, e.g. the new version FEniCSx [36].

### 3.1. PDE-type classes

The CUQIpy-FEniCS plugin defines PDE-type classes, see table 2, that represent PDE problems implemented using FEniCS. To view the underlying PDE-type class that is used in the 2D Poisson example, we call `print(A.pde)`, where `A` is the CUQIpy PDEModel defined in section 1.2, and obtain:

```
CUQI SteadyStateLinearFEniCSPDE.
PDE form expression:
def form(w,u,p):
    return ufl.exp(w)*ufl.inner(ufl.grad(u), ufl.grad(p))*ufl.dx
            - f*p*ufl.dx
```

Specifically, the Poisson PDE is represented by the `SteadyStateLinearFEniCSPDE` class. Similar to the core CUQIpy PDE-type classes, a CUQIpy-FEniCS PDE-type class contains a PDE `form`, which is a user-provided Python function that uses FEniCS syntax to express the PDE weak form; more discussion on building weak forms is provided in section 4 in the context of the EIT example. The Python function `form` inputs are `w`, the unknown parameter (log-conductivity in the 2D Poisson example), `u`, the state variable (or trial function), and `p`, the adjoint variable (or test function); and `f` is the FEniCS expression of the source term.

The CUQIpy-FEniCS PDE-type classes follow the interface defined by the core CUQIpy abstract PDE class by implementing the methods `assemble`, `solve`, and `observe`.

**Table 2.** Modules and classes of the CUQIpy-FEniCS plugin.

| Class name | Description |
| --- | --- |
| `cuqipy_fenics.pde` module: | |
| `FEniCSPDE` | A base (abstract) class that represents PDE problems defined using FEniCS |
| `SteadyStateLinearFEniCSPDE` | A class representation of steady state linear PDE problems defined using FEniCS |
| `cuqipy_fenics.geometry` module: | |
| `FEniCSContinuous` | A class representing FEniCS function spaces |
| `FEniCSMappedGeometry` | A class with additional mapping applied to the function values |
| `MaternKLExpansion` | A class that builds spectral representation of Matérn covariance operator on a given space, represented by a `FEniCSContinuous` geometry |
| `cuqipy_fenics.testproblem` module: | |
| `FEniCSDiffusion1D` | 1D diffusion PDE problem defined using FEniCS |
| `FEniCSPoisson2D` | 2D Poisson PDE problem defined using FEniCS |

The `assemble` method builds the discretized PDE system to be solved, from the provided PDE form. In the Poisson example, section 1.2, the system that results from discretizing the weak form of the PDE (1) using the FEM can be written as:

$$\boldsymbol{K}_w \boldsymbol{u} = \boldsymbol{f}, \tag{12}$$

where $\boldsymbol{K}_w$ is the discretized diffusion operator, given the discretized log-conductivity field $\boldsymbol{w}$. The vector $\boldsymbol{u}$ is the discretized PDE solution, the potential, and $\boldsymbol{f}$ is the discretized source term.

The `solve` method solves the linear system (12) using a FEniCS linear solver, that can be specified by the user. As discussed in section 2, this method represents the discretized solution operator $\boldsymbol{S}$, which in this case maps the log-conductivity $\boldsymbol{w}$ used in assembling $\boldsymbol{K}_w$ to the PDE solution $\boldsymbol{u}$.

Similarly, as discussed in section 2, the method `observe` represents the discretized observation operator $\boldsymbol{O}$. Since, in this case, the observations are obtained on the entire domain, $\boldsymbol{O}$ is just an identity operator that maps the full solution $\boldsymbol{u}$ to the observations $\boldsymbol{y} = \boldsymbol{u}$. In general, however, $\boldsymbol{O}$ can represent observing parts of the solution only, cf section 2, and/or a derived quantity of interest, cf section 4 for example.

The `SteadyStateLinearFEniCSPDE` class additionally implements the method `gradient_wrt_parameter` that computes the gradient of $\boldsymbol{O} \circ \boldsymbol{S}$ with respect to the parameter `w` in a given direction, using an adjoint-based approach [16]. The software design concept of the PDE `form` above and the adjoint-based gradient computation of the PDE `form` follows closely the approach used in hIPPYlib [42, 43].

For brevity we do not provide code for building the `SteadyStateLinearFEniCSPDE` object here, as it is provided by the CUQIpy-FEniCS test problem FEniCSPoisson2D, and stored in `A.pde`. How to build `PDE`-like objects is shown in section 2 for the core **CUQIpy**, and in section 4 for the **CUQIpy-FEniCS** plugin.

### 3.2. `Geometry`-type classes in CUQIpy-FEniCS

`Geometry`-type classes, as we discussed in section 2, mainly serve three purposes. First, they interface forward models with samplers and optimizers by providing a dual representation of variables: parameter value and function values representations. Second, they provide visualization capabilities for both representations. Lastly, they allow re-parameterizing the Bayesian inverse problem parameter, e.g. in terms of coefficients of expansion for a chosen basis. **CUQIpy-FEniCS** `Geometry`-type classes serve the same goals, see table 2 for a list of these classes.

There are two main data structures in **FEniCS**, `Function` and `Vector`. The former is a class representation of FEM approximation of continuous functions, and the latter is a class representation of the approximation coefficients of expansion. **CUQIpy-FEniCS** `Geometry`-type classes, subclassed from `cuqi.geometry.Geometry`, interpret these data structures and interface them with **CUQIpy**. Additionally, they provide plotting methods by seamlessly utilizing the **FEniCS** plotting capabilities. This enables **CUQIpy-FEniCS** to visualize function value representation of variables, as in figure 1(a), as well as, parameter value representation of variables, as in figure 1(f). The plotting implementation details are hidden from the user; and the user is provided with the **CUQIpy** simple plotting interface as shown for example in section 2.

**CUQIpy-FEniCS** `Geometry`-type classes provide useful parameterization and mapping functionalities. Here we discuss the `FEniCSContinuous` and the `MaternKLExpansion` `Geometry`-type classes, both are used in the 2D Poisson, the EIT, and the PAT examples. The most basic **CUQIpy-FEniCS** `Geometry`-type class is the `FEniCSContinuous` geometry which represents **FEniCS** FEM function spaces. We can write a FEM approximation $w^{\text{FEM}}(\boldsymbol{\xi})$ of a continuous function $w(\boldsymbol{\xi})$ as

$$w(\boldsymbol{\xi}) \approx w^{\text{FEM}}(\boldsymbol{\xi}) = \sum_{i=1}^{n_{\text{FEM}}} w_i^{\text{FEM}} e_i^{\text{FEM}}(\boldsymbol{\xi}), \tag{13}$$

where $\{e_i^{\text{FEM}}(\boldsymbol{\xi})\}_{i=1}^{n_{\text{FEM}}}$ are FEM basis functions defined on a given mesh, $\boldsymbol{w}^{\text{FEM}} = [w_1^{\text{FEM}}, w_2^{\text{FEM}}, \ldots, w_{n_{\text{FEM}}}^{\text{FEM}}]^{\mathsf{T}}$ is the vector of the corresponding FEM coefficients of expansions, and $n_{\text{FEM}}$ is the number of basis functions. The `FEniCSContinuous.par2fun` method converts a **NumPy** array to a **FEniCS** `Function` object representing $w^{\text{FEM}}(\boldsymbol{\xi})$. This is achieved by interpreting the array elements as the FEM expansion coefficients $\boldsymbol{w}^{\text{FEM}}$. The method `fun2par` converts a **FEniCS** `Function` objects representing $w^{\text{FEM}}(\boldsymbol{\xi})$ to a **NumPy** array of the FEM expansion coefficients $\boldsymbol{w}^{\text{FEM}}$. We denote by $\boldsymbol{G}_{\text{FEM}}$ the operator implemented by the `par2fun` method which maps $\boldsymbol{w}^{\text{FEM}}$ to $w^{\text{FEM}}(\boldsymbol{\xi})$. We use the FEM coefficient vector notation $\boldsymbol{w}^{\text{FEM}}$ when referring the FEM function $w^{\text{FEM}}(\boldsymbol{\xi})$, for simplicity. To create an object of the `FEniCSContinuous` class, which we use for example to represent the observations $\boldsymbol{y}$ in the Poisson example and refer to as `G_FEM`, we first define the **FEniCS** function space on which the parameter is represented

```
parameter_space = dl.FunctionSpace(mesh, "CG", 1)
```

mesh is the FEniCS computational mesh representing the physical domain of the problem, and parameter_space is a FEniCS first-order Lagrange polynomial space defined on mesh. We are now ready to create the FEniCSContinuous object as follows

```
G_FEM = FEniCSContinuous(parameter_space)
```

In some cases, re-parameterizing the Bayesian inverse problem parameter is needed to enforce certain type of solutions. One such re-parameterization, that is used in the Poisson example, is to enforce smooth solutions through a KL expansion. In CUQIpy-FEniCS, a KL parameterization can be represented by a MaternKLExpansion geometry. This geometry is used to approximate the FEM coefficient of expansion vector $w^{\text{FEM}}$ by a truncated KL expansion $w$:

$$w^{\text{FEM}} \approx w = \sum_{i=1}^{n_{\text{KL}}} x_i \sqrt{\lambda_i} e_i^{\text{KL}}. \tag{14}$$

Here, $x = [x_1, x_2, \ldots, x_{n_{\text{KL}}}]^\mathsf{T}$ is the KL-expansion coefficient vector, $\{\lambda_i\}_{i=1}^{n_{\text{KL}}}$ are a decreasing sequence of positive real numbers, and $\{e_i\}_{i=1}^{n_{\text{KL}}}$ is a set of FEM coefficient vectors of orthonormal functions. MaternKLExpansion constructs this KL expansion by discretizing a covariance operator—specifically, a Matérn-class covariance operator $(\frac{1}{\ell^2}I - \Delta)^{-(\frac{\nu}{2} + \frac{d}{4})}$ with $\ell > 0$, a smoothness parameter $\nu > 1$, and the physical domain spatial dimension $d = 1, 2$ or 3 [13]—on a FEM function space, parameter_space in this case. We then exploit FEniCS eigenvalue solvers to obtain the approximate eigenpairs $\{(\sqrt{\lambda_i}, e_i)\}_{i=1}^{n_{\text{KL}}}$. We refer to (14) as the KL parameterization of $w$ with the KL coefficients $x$. Note that choosing, $n_{\text{KL}} \ll n_{\text{FEM}}$ reduces the dimension of the parameter space which simplifies solving the Bayesian inverse problem and is typically an accurate approximation in representing smooth fields.

We denote by the operator $G_{\text{KL\_VEC}}$ the map from the KL expansion coefficients $x$ to the FEM expansion coefficients $w$. The MaternKLExpansion object thus represents the map $G_{\text{KL}} := G_{\text{FEM}} \circ G_{\text{KL\_VEC}}$. In the 2D Poisson example, section 1.2, the MaternKLExpansion is internally created by the FEniCSPoisson2D test problem as

```
G_KL = MaternKLExpansion(G_FEM, length_scale=0.1, num_terms=32)
```

and it is used as the domain geometry of the model A to approximately parametrize the log-conductivity $w^{\text{FEM}}(\xi)$ by KL expansion coefficients $x$. The MaternKLExpansion class additionally implements the method gradient which computes the gradient of the map $G_{\text{KL}}$ with respect to the coefficients $x$ in a given direction.

### 3.3. Integration into CUQIpy through the PDEModel class

The CUQIpy-FEniCS PDE-type and Geometry-type objects provide the building blocks required to create the forward map $A$, e.g. (2). The CUQIpy PDEModel combines these FEniCS-dependent objects and interface them to the core CUQIpy library. We run print(A), where A is the CUQIpy model defined in section 1.2 to see its contents:

```
CUQI PDEModel: MaternKLExpansion on FEniCSContinuous(1089,) ->
               FEniCSContinuous(1089,).
    Forward parameters: ['x'].
    PDE: SteadyStateLinearFEniCSPDE.
```

We note its domain geometry (in the first line) corresponds to `G_KL` and its range geometry (second line) is `G_FEM`. We write $A$ in terms of its components as:

$$A(x) = (O \circ S \circ G_{\text{KL}})(x). \tag{15}$$

`PDEModel` provides the `forward` method that corresponds to applying $A$. Additionally, it provides a `gradient` method to compute the gradient of $A$ with respect to the parameter $x$, if the underlying `Geometry` and `PDE` objects have gradient implementation. This enables gradient-based posterior sampling such as by the NUTS sampler which we use in section 1.2.

The following section gives an elaborate case study of using CUQIpy-FEniCS to solve an EIT problem. We provide details of constructing the software components needed to build and solve the problem using CUQIpy equipped with CUQIpy-FEniCS.

## 4. CUQIpy-FEniCS example: EIT

EIT is an imaging technique of inferring the conductivity of an object from measurements of the electrical current on the boundary. This is a non-invasive approach in medical imaging for detecting abnormal tissues. Similar techniques are used in many other applications, such as industrial inspection. The underlying mathematical model for EIT is an elliptic PDE. Such PDEs are some of the most popular models for PDE-based inverse problems, e.g. in modeling subsurface flow in a porous medium, and inverse steady-state heat transfer problem with unknown heat conductivity. Hence, the EIT model can be modified for use in a wide range of PDE-based inverse problems.

Inferring discontinuous fields is a common problem in many inverse problems applications. Such fields are used, for example, to model tumors in medical imaging applications, abnormalities in fault detection applications, and inhomogeneity in geophysics applications. A classic approach to incorporate such fields into a Bayesian inverse problem is to use a Markov random field (MRF)-type prior, discussed in detail in [35]. However, such priors often result in a large set of parameters, yielding inefficient numerical uncertainty quantification methods for fine discretization levels. Recently, a new class of Bayesian priors has emerged where a discontinuous field is constructed using a non-linear transformation of a continuous prior [2, 13, 24]. A continuous prior can be constructed with a relatively few parameters, e.g. using a KL expansion. A deterministic non-linear transformation is then chosen prior to inference to capture the properties of the unknown field being inferred. Therefore, efficient numerical inference methods can be constructed. The level-set prior [24] is one such prior where discontinuities in the field are defined to be the zeros level-set of a smooth Gaussian random field, e.g. a KL expansion. In this EIT example, we utilize the Bayesian level-set approach to perform uncertainty quantification for the EIT problem.

*4.1. Mathematical model of EIT*

We follow the EIT model presented in [37] in 2D. Let $\Gamma \subset \mathbb{R}^2$ be the open unit disk representing an object, and let $\sigma : \bar{\Gamma} \to \mathbb{R}$ represent the conductivity of the object, where $\bar{\Gamma}$ is the closure of the set $\Gamma$. Suppose we impose a sequence of electric potentials $g_k(\boldsymbol{\xi})$, $k \in \mathbb{N}$, at the boundary $\partial \Gamma$. We can then find the distribution of the electric potential $u_k$, associated to $g_k$, inside the object from the elliptic PDE

$$-\nabla \cdot (\sigma(\boldsymbol{\xi}) \nabla u_k(\boldsymbol{\xi})) = 0, \qquad \boldsymbol{\xi} \in \Gamma, k \in \mathbb{N}^+ \quad (16a)$$

$$u_k(\boldsymbol{\xi}) = g_k(\boldsymbol{\xi}) = \sin\left(k \arctan\left(\frac{\xi^2}{\xi^1}\right)\right), \qquad \boldsymbol{\xi} \in \partial\Gamma, k \in \mathbb{N}^+. \quad (16b)$$

Here, $\sigma(\boldsymbol{\xi})$ is the conductivity as a function of the Cartesian coordinates $\boldsymbol{\xi} = [\xi^1, \xi^2]^\mathsf{T}$, and $k \in \mathbb{N}^+$ is some spatial frequency of the boundary electric potential. Note that the boundary condition (16*b*) for the electric potential is chosen from [5] which is an approximation of the full-electrode model introduced in [37].

The EIT problem is the inverse problem of inferring the conductivity $\sigma$ by imposing multiple boundary potentials $g_k(\boldsymbol{\xi})$, e.g. for $k = 1, 2, 3, \ldots$, and measuring the corresponding current $y_k(\boldsymbol{\xi})$, defined by

$$y_k(\boldsymbol{\xi}) := \frac{\sigma(\boldsymbol{\xi}) \partial u_k(\boldsymbol{\xi})}{\partial \mathbf{n}}, \quad (17)$$

at the boundary $\partial \Gamma$. Here, $\mathbf{n}$ is the outward unit vector, orthogonal to the boundary $\partial \Gamma$. This EIT model corresponds to the Dirichlet-to-Neumann EIT model, also known as the shunt model [37].

We are interested in piece-wise constant conductivity $\sigma$ with background conductivity value $\sigma^- = 1$ and foreground conductivity value $\sigma^+ = 10$. This contrast between foreground and background is a common difference between a healthy and an unhealthy (e.g. cancerous) tissue [5]. We also assume that the foreground represents an inclusion far from the boundaries, simplifying the boundary measurement to

$$y_k(\boldsymbol{\xi}) = \frac{\partial u_k(\boldsymbol{\xi})}{\partial \mathbf{n}}, \quad \boldsymbol{\xi} \in \partial\Gamma, \quad (18)$$

since $\sigma(\boldsymbol{\xi}) \equiv \sigma^- = 1$ on the boundary. We define the parameter-to-solution operator $\mathcal{S}_k$ as the mapping from the conductivity $\sigma(\boldsymbol{\xi})$ to the solution $u_k(\boldsymbol{\xi})$ of the PDE (16). We also define the observation operator $\mathcal{O}$ that maps the PDE solution $u_k(\boldsymbol{\xi})$ to the boundary current measurement $y_k(\boldsymbol{\xi})$, with $\boldsymbol{\xi} \in \partial\Gamma$. Note that the observation operator $\mathcal{O}$ does not explicitly depend on the frequency $k$.

In practice only a finite number of frequencies $k$, in (16), is considered. In this section we only consider $k = 1, 2, 3$ and 4.

*4.2. Finite element discretization and FEniCS implementation of EIT*

Let $H^1(\Gamma)$ [1] be the Hilbert space in which we expect the solution $u_k$ of equation (16) to belong. We now reformulate (16) to obtain an elliptic PDE with homogeneous Dirichlet boundary conditions. This can be achieved e.g. using the lifting method [32]. This approach simplifies the finite element approximation of (16).

We define lifting functions $u_k^{\text{lift}} \in H^1(\Gamma)$, $k = 1, 2, 3$ and 4, that satisfy the boundary input in (16), i.e. $u_k^{\text{lift}}(\boldsymbol{\xi})|_{\partial\Gamma} = u_k(\boldsymbol{\xi})|_{\partial\Gamma} = g_k(\boldsymbol{\xi})$, and vanishes away from the boundary. Introducing a new variable $v_k = u_k - u_k^{\text{lift}}$ allows us to reformulate (16) as

$$-\nabla \cdot (\sigma(\boldsymbol{\xi}) \nabla v_k(\boldsymbol{\xi})) = \Delta u_k^{\text{lift}}(\boldsymbol{\xi}), \qquad \boldsymbol{\xi} \in \Gamma, \ k = 1, 2, 3, 4, \qquad (19a)$$
$$v_k(\boldsymbol{\xi}) = 0, \qquad \boldsymbol{\xi} \in \partial\Gamma. \qquad (19b)$$

The potential $u_k$ is now recovered from the relation $u_k = v_k + u_k^{\text{lift}}$, for $k = 1, 2, 3$ and 4.

Taking test function $t(\boldsymbol{\xi}) \in H^1(\Gamma)$, we form the weak formulation [32] of (19) as

$$\int_\Gamma \sigma \nabla v_k \cdot \nabla t \, d\boldsymbol{\xi} = -\int_\Gamma \nabla u_k^{\text{lift}} \cdot \nabla t \, d\boldsymbol{\xi}. \qquad (20)$$

Similarly, we let $H_{\text{p}}^1(\partial\Gamma)$ be the space which the observation function $y_k$ belong to. Here, the subscript 'p' denotes the Hilbert space of periodic functions. Taking a test function $w \in H_{\text{p}}^1(\partial\Gamma)$, we can form the weak form for the boundary measurement as

$$\int_{\partial\Gamma} y_k w \, d\boldsymbol{\xi} = \int_{\partial\Gamma} \frac{\partial v_k + \partial u_k^{\text{lift}}}{\partial \mathbf{n}} w \, d\boldsymbol{\xi}. \qquad (21)$$

Here, we used the relation $\partial u_k / \partial \mathbf{n} = \partial(v_k + u_k^{\text{lift}}) / \partial \mathbf{n}$. Note that due to the lifting approach, the observation now depends on the frequency $k$. We emphasize this by introducing the subscript $k$ for the observation operator, i.e. we define $\mathcal{O}_k$ to be the mapping from $v_k$ to $y_k$.

We discretize the domain $\Gamma$ using a triangulated mesh. Furthermore, we choose first-order Lagrangian polynomial functions [32] to approximate the basis functions of $H^1(\Gamma)$ and $H_{\text{p}}^1(\partial\Gamma)$. We implement the left-hand-side and the right-hand-side of (20) using **FEniCS** as

```
form_lhs = lambda sigma, v, t:
        dl.inner(sigma*dl.grad(v), dl.grad(t))*dl.dx
form_rhs1 = lambda sigma, t:
        -dl.inner(dl.grad(u_lift_1), dl.grad(t))*dl.dx
```

Here, the functions `form_lhs` and `form_rhs1` return the **FEniCS** weak forms of the left-hand side and the right-hand side (for $k = 1$) of (20), respectively. Furthermore, the **FEniCS** function `u_lift_1` contains the user-defined lifting function. We refer the reader to the codes accompanying this paper for more details. Note that since $v_k$ is the solution to the PDE (20), we may use the same `form_lhs` for all frequencies $k = 1, 2, 3$ and 4. We construct similar functions `form_rhs2`, `form_rhs3`, and `form_rhs4` for the input frequencies $k = 2, 3$ and 4.

We now implement the observation function (21). Let `give_bnd_vals` be a Python function that computes function values at the boundaries of $\Gamma$. The observation function then takes the form

```
def observation1(sigma, v1):
    obs_form = dl.inner(dl.grad(v1 + u_lift_1), n)*w*dl.ds
    assembled_form = dl.assemble(obs_form)
    boundary_values = give_bnd_vals(assembled_form)
    return boundary_values
```

Here, `n` is a **FEniCS** vector containing the unit outward normal vectors to the cell boundaries and `v1` is a **FEniCS** function of the solution $v_1$, and `w` is a **FEniCS** test function. We construct similar functions `observation2`, `observation3`, and `observation4` for the input frequencies $k = 2, 3$ and $4$.

The FEM discretization of (20) results in a finite-dimensional system of equations

$$\boldsymbol{M}_{\boldsymbol{\sigma}} \boldsymbol{v}_k = \boldsymbol{b}_k, \qquad k = 1, 2, 3, 4, \tag{22}$$

and the discretized observation model

$$\boldsymbol{y}_k = \boldsymbol{O}_k(\boldsymbol{v}_k), \qquad k = 1, 2, 3, 4. \tag{23}$$

Here, $\boldsymbol{\sigma}$ and $\boldsymbol{v}_k$ are the FEM expansion coefficients for $\sigma$ and $v_k$, respectively. Furthermore, $\boldsymbol{M}_{\boldsymbol{\sigma}}$ is the FEM mass matrix, i.e. the discretization of the elliptic operator $\nabla \cdot \sigma \nabla$, that depends on the unknown parameter $\boldsymbol{\sigma}$, and $\boldsymbol{b}_k$ is the right-hand-side vector containing the estimated integrals in the right-hand-side of (20). Furthermore, $\boldsymbol{y}_k$ is the observation vector, and $\boldsymbol{O}_k$ is a discretization of the observation map $\mathcal{O}_k$.

We now define the EIT forward maps $\boldsymbol{A}_k$ to be

$$\boldsymbol{y}_k = \boldsymbol{A}_k(\boldsymbol{\sigma}) := \boldsymbol{O}_k \circ \boldsymbol{S}_k(\boldsymbol{\sigma}), \qquad k = 1, 2, 3, 4, \tag{24}$$

where $\boldsymbol{S}_k$ is the FEM discretization of $\mathcal{S}_k$ discussed above.

### 4.3. Parameterization of the conductivity $\sigma$

In this section we consider the level-set parameterization for the conductivity $\sigma$ proposed in [13]. This approach comprises multiple layers of parameterization. In this section we use the `geometry` module in CUQIpy-FEniCS to implement such layered parameterization.

Let us first define the FEniCS function space in which we expect $\boldsymbol{\sigma}$ to belong

```
parameter_space = dl.FunctionSpace(mesh, "CG", 1)
```

Here, `mesh` is the computational mesh used for FEniCS. Recall that `parameter_space` is a FEniCS function space with linear continuous elements. Similarly, we can define a FEniCS function space `solution_space` on which the solutions $v_k$ defines a function.

Now we define a discrete Gaussian random field $\boldsymbol{r}$ defined on $\Gamma$, i.e. realizations of $\boldsymbol{r}$ are FEM expansion coefficients of a random functions defined on $\Gamma$. One way to define such a random function is to use a truncated KL-expansion with a Matérn covariance, as discussed in section 3.2, to approximate $\boldsymbol{r}$ the same way $\boldsymbol{w}$ is approximated in (14).

To construct the geometry associated to the KL parameterization, we first consider the operator $\boldsymbol{G}_{\text{FEM}}$ to be the map from FEM expansion coefficients to a FEniCS function (see section 3.2). The corresponding geometry is defined as

```
G_FEM = FEniCSContinuous(parameter_space)
```

We now construct the geometry associated with the KL parametrization. This geometry is associated with the operator $\boldsymbol{G}_{\text{KL}}$ defined in section 3.2.

```
G_KL = MaternKLExpansion(G_FEM, length_scale=0.2, num_terms=64)
```

Here `length_scale` is the length scale constant of the Matérn covariance and `num_terms` is $n_{\text{KL}}$, the number of terms in the KL expansion. The geometry `G_KL` is now the implementation of $\boldsymbol{G}_{\text{KL}}$ which maps $\boldsymbol{x}$, the vector containing the KL expansion coefficients in (14) to the vector $\boldsymbol{r}$. Note that we used `parameter_space` as the FEniCS function space associated with $\boldsymbol{r}$.

Now to relate the Gaussian function $r$ to the conductivity $\sigma$ we define the Heaviside function [45] $G_{\text{Heavi}}$ as an additional layer of parameterization

$$\sigma = G_{\text{Heavi}}(r) := \frac{1}{2}\left(\sigma^+(1-\text{sign}(r)) + \sigma^-(1+\text{sign}(r))\right). \tag{25}$$

This map $G_{\text{Heavi}}$ constructs a piece-wise constant conductivity $\sigma$. Note that the Heaviside map must be applied to the function $r$, the FEM function associated to $r$, which constructs a conductivity $\sigma$. However, in the case of linear Lagrangian FEM elements, we can directly apply the Heaviside map to the expansion coefficients $r$ and obtain $\sigma$ as in (25).

We can construct this parameterization in CUQIpy-FEniCS as

```
G_Heavi = FEniCSMappedGeometry(G_KL, map=heaviside)
```

Here, `heaviside` is a Python function that applies the Heaviside map (25), see the companion code for more implementation details. By passing `map = heaviside`, `FEniCSMappedGeometry` applies `heviside` to `G_KL`.

We redefine the forward operators to use the parameterizations discussed:

$$A_k = O_k \circ S_k \circ G_{\text{Heavi}} \circ G_{\text{KL}}, \qquad k = 1, 2, 3, 4. \tag{26}$$

We define the range geometry as a `Continuous1D` geometry:

```
G_cont = Continuous1D(m)
```

where `m` is the dimension of any observation vector $y_1, y_2, y_3$ or $y_4$. In the experiments in this section we set `m = 94`.

### 4.4. `PDEmodel` for the EIT problem

Now we have all the components to define a CUQIpy-FEniCS PDE object. We first create a PDE form by combining the left-hand-side and right-hand-side forms, defined in section 4.2, in a Python tuple as

```
PDE_form1 = (form_lhs, form_rhs1)
```

Since (19) is a steady state linear PDE, we use the `SteadyStateLinearFEniCSPDE` class to define this PDE in CUQIpy-FEniCS.

```
PDE1 = SteadyStateLinearFEniCSPDE(
        PDE_form1, mesh, solution_space,
        parameter_space, zero_bc,
        observation_operator=observation1,
        reuse_assembled=True)
```

Recall that `solution_space` is the FEniCS space for the solution $v_k$ in (22), `zero_bc` is the FEniCS implementation of the homogeneous Dirichlet boundary conditions for (19), and `observation1` is the observation Python function defined in section 4.2. The key argument `reuse_assembled = True` informs CUQIpy-FEniCS to store and reuse matrix factors of $M_\sigma$, for a particular $\sigma$, when solving the system (22). This provides a significant computational acceleration.

We note that PDE problems (19) for frequencies $k = 2,3$ and 4 differ from the frequency $k = 1$ only in the right-hand-side term and the observation operator. We can exploit this to construct PDE2, for frequency $k = 2$, as

```
PDE2 = PDE1.with_updated_rhs(form_rhs2)
PDE2.observation_operator = observation2
```

And similarly we construct PDE3 and PDE4. Note that the matrix factorization of $M_\sigma$ is shared among PDE1, PDE2, PDE3, and PDE4.

Now we can create a `PDEModel` that represents the forward operator $A_1$ and includes information about the parameterization of $\sigma$.

```
A1 = PDEModel(PDE1, range_geometry=G_cont, domain_geometry=G_Heavi)
```

Similarly we define A2, A3, and A4 for the input frequencies $k = 2,3$ and 4 corresponding to the forward operators defined in (26).

### 4.5. Bayesian formulation and solution

In this section we formulate the EIT problem in a Bayesian framework. Let $x$ be a vector containing the expansion coefficients $\{x_i\}_{i=1}^{n_{\text{KL}}}$ in (14). In the Bayesian formulation of the EIT problem the posterior distribution of the unknown parameter $x$ is the conditional probability distribution of $x$ given observed data $y_1^{\text{obs}}, y_2^{\text{obs}}, y_3^{\text{obs}}$ and $y_4^{\text{obs}}$. Here, we assume white Gaussian data noise. This Bayesian problem takes the form

$$x \sim \text{Gaussian}\left(\mathbf{0}, I_{n_{\text{KL}}}\right),$$
$$y_k \sim \text{Gaussian}\left(A_k(x), s_{\text{noise}}^2 I_m\right), \qquad k = 1,2,3,4.$$

Here, $s_{\text{noise}}$ is the standard deviation of the data distribution. We use CUQIpy to implement these distributions.

```
x = Gaussian(np.zeros(n_KL), 1, geometry=G_Heavi)
y1 = Gaussian(A1(x), s_noise**2, geometry=G_cont)
```

and similarly we define y2, y3, and y4 for $k = 2,3$ and 4. We pass the argument `geometry = G_Heavi` when initializing x to specify that samples from this distribution follow the parameterization discussed in section 4.3. We can sample from the prior distribution and plot the samples using the following:

```
prior_sample = x.sample(5)
prior_sample.plot()
```

Examples of prior samples can be found in figure 7. Note that CUQIpy-FEniCS visualizes these samples as FEniCS functions.

To create simulated data for this EIT problem, we consider the conductivity field $\sigma^{\text{true}}$ comprising 3 circular inclusions. The coordinates of the centers of the inclusions are $(0.5, 0.5)$, $(-0.5, 0.6)$ and $(-0.3, -0.3)$ with radii 0.2, 0.1 and 0.3, respectively. We also assume conductivity values of $\sigma^+ = 10$ and $\sigma^- = 1$ inside and outside of inclusions, respectively. We can obtain the FEM expansion coefficients $\sigma^{\text{true}}$ by projecting $\sigma^{\text{true}}$ onto the FEM basis. Note that we introduce an approximation in this projection. The true and the projected conductivity phantoms are presented in figure 6(a). Note that this conductivity field is not sampled from the prior distribution, and thus, there is no true parameters $\boldsymbol{x}^{\text{true}}$ that gives rise to the exact conductivity $\boldsymbol{\sigma}^{\text{true}}$.

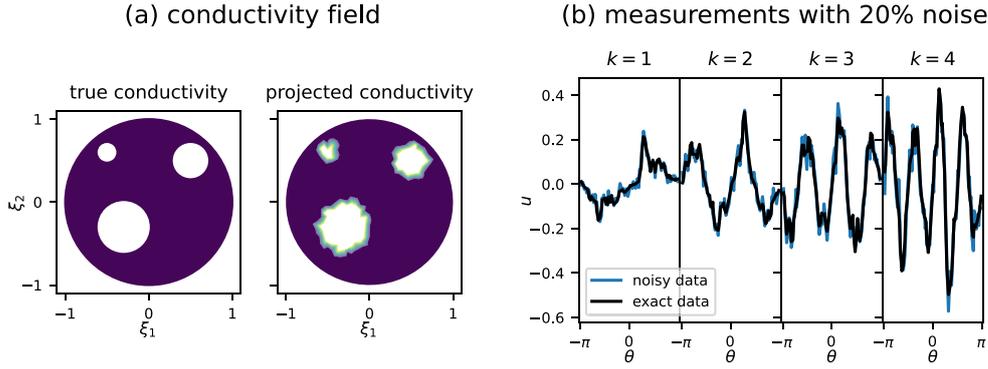

**Figure 6.** (a) The true (but assumed unknown) conductivity field $\sigma$ (left) and the projected conductivity field onto the FEM space (right). (b) The exact boundary measurement values and the noisy boundary measurement values for frequencies $k = 1,2,3,4$. Here we only present the data collected with 20% noise level.

Data is created by adding additive Gaussian noise, with standard deviation $s_{\text{noise}}$, to $\boldsymbol{y}_k^{\text{exact}} := A_k(\boldsymbol{\sigma}^{\text{true}})$, for $k = 1,2,3,4$, at noise levels 5%, 10% and 20%. The true and the noisy data with 20% noise level are presented in figure 6(b).

Now we have all the components we need to create a posterior distribution. We first define the joint distribution

$$p(\boldsymbol{x},\boldsymbol{y}_1,\boldsymbol{y}_2,\boldsymbol{y}_3,\boldsymbol{y}_4) = p(\boldsymbol{x})\,p(\boldsymbol{y}_1|\boldsymbol{x})\,p(\boldsymbol{y}_2|\boldsymbol{x})\,p(\boldsymbol{y}_3|\boldsymbol{x})\,p(\boldsymbol{y}_4|\boldsymbol{x}), \qquad (27)$$

where $p(\boldsymbol{x})$ is the prior probability density function (PDF) and $p(\boldsymbol{y}_k|\boldsymbol{x})$, for $k = 1,2,3$ and 4, are the data distribution PDF. We obtain the posterior distribution by conditioning the joint distribution on the data $\boldsymbol{y}_k^{\text{obs}}$, for $k = 1,2,3$ and 4. In CUQIpy this translates to

```
joint = JointDistribution(x, y1, y2, y3, y4)
posterior = joint(y1=y1_obs, y2=y2_obs, y3=y3_obs, y4=y4_obs)
```

In this test case, we use the standard Metropolis-Hastings (MH) algorithm [35, § 2] to sample from the posterior. We pass the posterior as an argument in the initialization and then compute $10^6$ samples using this sampler.

```
sampler = MH(posterior)
posterior_samples = sampler.sample_adapt(1000000)
```

In what remains in this section we discuss how to use `posterior_samples` in CUQIpy-FEniCS to visualize the posterior distribution.

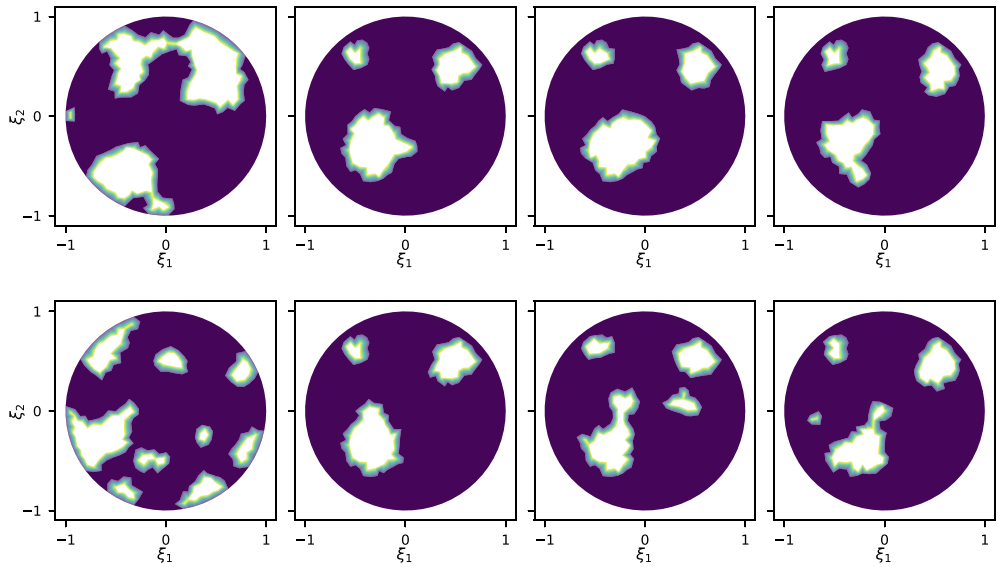

**Figure 7.** Samples from the prior and posterior distributions of $\sigma$. First column: prior samples. Second, third and fourth columns: posterior samples of $\sigma$ with 5%, 10% and 20% noise level, respectively.

*4.6. Posterior samples, post-processing and visualization*

We analyze and visualize the posterior using CUQIpy equipped with CUQIpy-FEniCS geometries. We first plot some of the posterior samples using the `plot` method.

```
posterior_samples.plot()
```

This command chooses 5 random posterior samples and plots them. We provide these posterior samples of the conductivity field $\sigma$ in figure 7 (second to fourth column), for noise levels 5%, 10% and 20%, respectively, only two samples are shown for each case for brevity. We see that the reconstruction in the samples degrades, compared to the true conductivity $\sigma^{\text{true}}$, for higher noise levels. In addition, we see the discrepancy between the samples and $\sigma^{\text{true}}$ happens near the center of the domain.

Now we want to estimate and visualize the posterior mean as an estimate for the conductivity field $\sigma$. We can achieve this using the command

```
posterior_samples.plot_mean()
```

Note that `posterior_samples` is equipped with the `G_heavi` geometry. Therefore, `plot_mean` will apply this geometry, i.e. the parameterization $\boldsymbol{G}_{\text{Heavi}} \circ \boldsymbol{G}_{\text{KL}}$ to the posterior mean. The mean conductivity field is provided in figure 8(a). We see that for increased noise level, the posterior mean less resembles the true conductivity field $\sigma^{\text{true}}$.

We use the point-wise variance of the posterior samples as a method for quantifying the uncertainty in the posterior. We can achieve this in **CUQIpy-FEniCS** by

```
posterior_samples.funvals.vector.plot_variance()
```

Here, the `Samples` property `funvals` converts the parameter samples to function value samples and return them in a new `Samples` object, i.e. it applies the map $G_{\text{Heavi}} \circ G_{\text{KL}}$ to generate the function value samples. Similarly, the `Samples` property `vector` converts these samples to a vector representation, the DOF of the **FEniCS** functions in this case. Eventually, the variance is computed over this vector representation and then plotted as a **FEniCS** function. The point-wise variance is presented in figure 8(b). We see that the uncertainty in the reconstruction is associated with the boundaries of the inclusions, as well as, distance from the domain boundary. Furthermore, adding noise increases the level of uncertainty. This is consistent with findings of [13].

Finally we visualize the posterior for the expansion coefficients $x$. We use `plot_ci` method to visualize the posterior mean and the 95% CIs associated with the parameters. To indicate that we are visualizing the posterior for the coefficient (parameter) $x$, we pass the argument `plot_par = True` to the `plot_ci` method.

```
posterior_samples.plot_ci(95, plot_par=True)
```

In figures 9(a)–(c) we present the CI plots for noise levels 5%, 10% and 20%, respectively. We see that for the case with 20% noise level, the mean of $x_i$ is close to zero, for larger indices $i$. This suggests that for higher noise levels, the value of $x_i$ follows the prior distribution and the information associated to these coefficients is lost in the data distribution. This is not the case for smaller noise levels, e.g. 5%.

## 5. PAT through user-defined PDE models

In many applications of UQ for inverse problems a well-developed forward problem solver is created by the user. Therefore, it is of high interest that **CUQIpy** and **CUQIpy-FEniCS** can incorporate such black-box forward solvers.

In this section we discuss how to use a black-box code in **CUQIpy** and **CUQIpy-FEniCS** to quantify uncertainties for inverse problems with PDEs. In addition, we discuss how to exploit geometries in **CUQIpy-FEniCS**, in such cases, to visualize uncertainties, without modifying the original black-box software.

To demonstrate the user-defined features of **CUQIpy** and **CUQIpy-FEniCS**, we consider a 1D PAT problem [41]. In such problems, a short light pulse is illuminated onto an object to create a local initial pressure distribution. This pressure distribution then propagates in the object in the form of ultrasound waves. The PAT problem is then to reconstruct the initial pressure distribution from time-varying ultrasound measurements. For the 1D variant, we consider a 1D pressure profile with $r = 2$ ultrasound sensors to measure pressure variations.

### 5.1. Mathematical model of PAT

Let us consider an infinitely long 1D acoustic object with homogeneous acoustic properties (homogeneous wave speed). Assuming that the illumination duration via a light pulse is negligible compared to the speed of wave propagation, we can approximate the propagation of waves in the object by the hyperbolic PDE (linear wave equation)

$$\frac{\partial^2 u(\tau,\xi)}{\partial \tau^2} = \frac{\partial^2 u(\tau,\xi)}{\partial \xi^2}, \qquad 0 < \tau \leqslant 1,\, \xi \in \mathbb{R}, \tag{28a}$$

$$u(0,\xi) = g(\xi), \qquad \xi \in \mathbb{R}, \tag{28b}$$

$$\frac{\partial u(0,\xi)}{\partial \tau} = 0, \qquad \xi \in \mathbb{R}. \tag{28c}$$

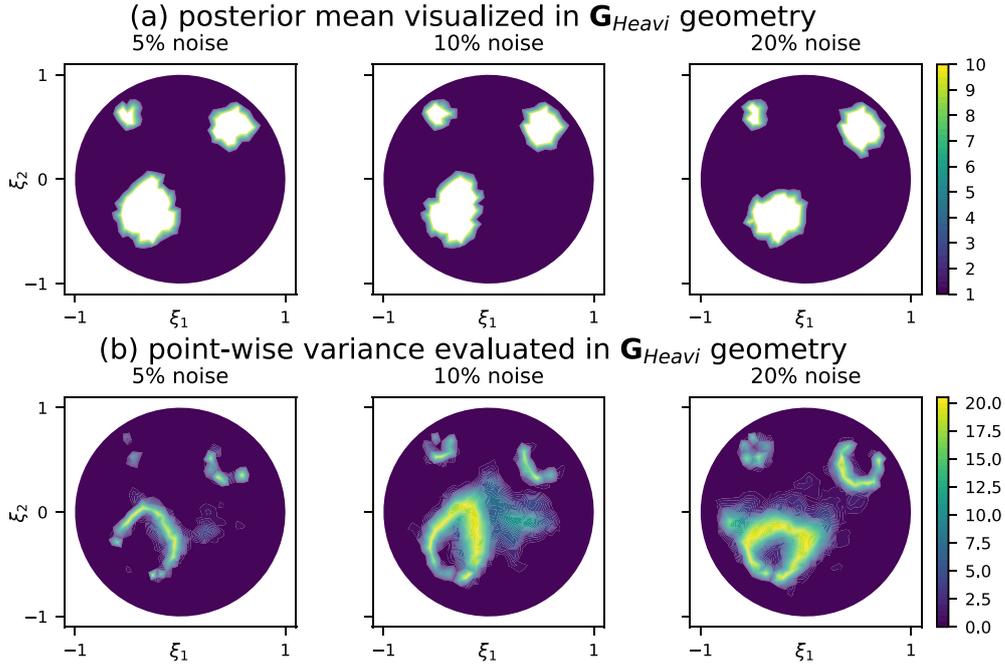

**Figure 8.** Estimated conductivity field $\sigma$ with uncertainty estimates, visualized as Heaviside-mapped KL expansion. (a) Posterior mean (b) point-wise variance.

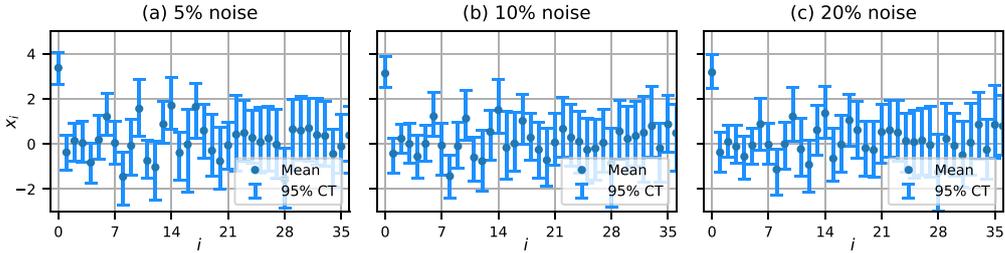

**Figure 9.** Estimation of $x$, i.e. the first 35 coefficients of the KL expansion in (14) and their uncertainty. (a) 5% noise (b) 10 % noise (c) 20 % noise.

Here, $u$ is the pressure distribution, $g$ the initial pressure distribution, and $\tau$ the time.

We place 2 pressure sensors at locations $\xi_L = 0$ and $\xi_R = 1$ that record the pressure over time. We define continuous measurements

$$y_L(\tau) := u(\xi_L, \tau), \qquad y_R(\tau) := u(\xi_R, \tau), \qquad 0 < \tau \leqslant 1. \tag{29}$$

Let $y := [y_L, y_R]$ and $r = 2$ in case of the availability of full-data and $y := y_L$ and $r = 1$ in case only partial-data is available. The inverse problem is to find the initial pressure $g$ from measurements $y$.

We define the parameter-to-solution operator $\mathcal{S}$ of the PAT problem to be the map that maps $g(\xi)$ to the solution $u(\tau, \xi)$ for all $0 < \tau \leqslant 1$ and $\xi \in \mathbb{R}$. Furthermore, we define the observation operator $\mathcal{O}$ to be the mapping that maps (slices) the solution $u(\tau, \xi)$ to the measurements $y$. Note that the forward operator

$$\mathcal{A} := \mathcal{O} \circ \mathcal{S}, \tag{30}$$

is a linear operator. We are interested in recovering the initial pressure profile $g(\xi)$ in the domain $\Gamma = [\xi_L, \xi_R] = [0, 1]$ from the full data, where pressure measurements from both $\xi_L$ and $\xi_R$ is available. In addition, we also investigate reconstructing $g$ from partial data, where pressure measurement is only available at $\xi_L$.

### 5.2. CUQIpy implementation of PAT

In this section we assume that a discretization $S$ of $\mathcal{S}$ is available that approximately solves the wave equation (28). Furthermore, we assume that the input of $S$ is a vector $g$ which represents a discretization of $g$. Note that the discretization $S$ is a generic discretization and we treat it as a black-box solver. However, the particular discretization used in this paper is available in the accompanying codes.

We consider a discretization $O$ of the observation operator $\mathcal{O}$. We choose a measurement frequency $f = 250$ and construct discrete measurement vector $y$ comprising snapshots $[u(\xi_L, i/f), u(\xi_R, i/f)]$ and $u(\xi_L, i/f)$ for the full-data and partial-data measurement, respectively, where $i = 1, \ldots, m$ and $m = 250$.

We consider a simple Bayesian problem where we assume the components of $g$ have a Gaussian distribution and the pressure measurements are corrupted by additive white Gaussian noise. We can formulate this problem as

$$g \sim \text{Gaussian}\left(\mathbf{0}, \mathbf{I}_{n_g}\right),$$
$$y \sim \text{Gaussian}\left(A(g), s_{\text{noise}}^2 \mathbf{I}_{rm}\right).$$

Here, $n_g = 121$ is the size of $g$, and $s_{\text{noise}} = 0.125$ is the standard deviation of the noise. Furthermore, $A$ represents the user's discretization of the forward operator $\mathcal{A}$, and `PAT` is a black-box Python function applying the photo-acoustic forward model $A$.

We set up the Bayesian Problem with Gaussian prior and data distribution:

```
g = Gaussian(np.zeros(n_g), 1)
y = Gaussian(PAT(x), s_noise**2)
```

Note that CUQIpy, by default, considers Continuous1D geometry for the initial pressure. Therefore, we do not explicitly define it. When incorporating black-box forward solvers within CUQIpy, the user has access to complex priors for g, following examples in the companion paper [35]. As an example, when the location of jumps in the initial pressure is the quantity of interest in the inverse problem, we can consider a Markov-random-field-type prior. Using such prior distributions is discussed in [35].

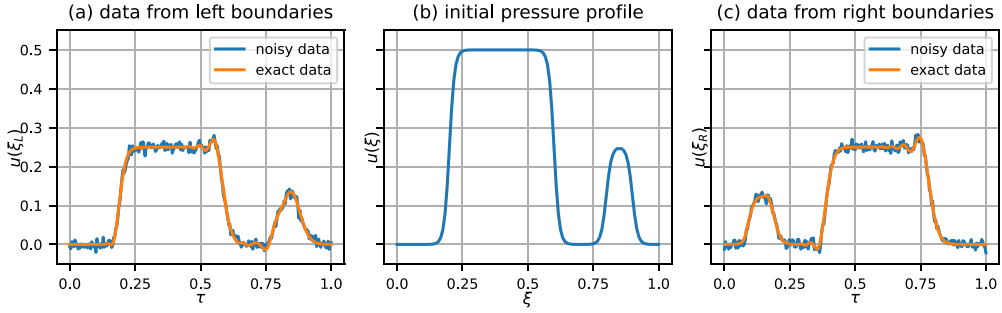

**Figure 10.** (b) True initial pressure for the PAT problem. (a) and (c) noisy and noise-free data collected for the PAT problem with sensors at $\xi_L = 0$ and $\xi_R = 1$, respectively.

We consider now a known true initial pressure $g^{\text{true}}$ and its discretization array `g_true` from which we can construct noisy measurements $\boldsymbol{y}^{\text{obs}}$

```
y_obs = y(g=g_true)
```

Figure 10 shows the initial pressure distribution and exact and noisy pressure measurements for sensors located at $\xi_L = 0$ and $\xi_R = 1$.

Instead of constructing the posterior and sampling it, we wish to demonstrate how to formulate this same problem using the CUQIpy-FEniCS plugin. Since sampling is done in the same way, we demonstrate it at the end of the following section.

### 5.3. CUQIpy-FEniCS implementation of PAT

Here, we assume that `PAT` is a Python function with a FEniCS implementation of the forward operator $\boldsymbol{A}$. We also assume this function take a FEniCS function g as input and computes the boundary pressure measurement y. Let us first define the FEniCS function space where $\boldsymbol{g}$ defines a function

```
parameter_space = dl.FunctionSpace(mesh, "CG", 1)
```

Here, `mesh` is a discretization of the real line and `parameter_space` is a FEniCS function space with first order Lagrangian hat-functions. We parameterize $\boldsymbol{g}$ with a KL-expansion with a Matérn covariance associated to the map $\boldsymbol{G}_{\text{KL}}$ (see sections 3.2 and 4.3). We now redefine the forward operator

$$\boldsymbol{A} = \boldsymbol{O} \circ \boldsymbol{S} \circ \boldsymbol{G}_{\text{KL}}. \tag{31}$$

We set up a geometry for the KL-expansion with CUQIpy-FEniCS (see section 3.2) as

```
G_FEM = FEniCSContinuous(parameter_space)
G_KL = MaternKLExpansion(G_FEM,
        length_scale=0.1, nu=0.75, num_terms=100)
```

The geometry `G_KL` is the implementation of $\boldsymbol{G}_{\mathrm{KL}}$ with Matérn length scale constant $\ell = 0.1$ and regularity constant $\nu = 0.75$, and $n_{\mathrm{KL}} = 100$ terms. We set up a `FEniCSMappedGeometry` with a map function `prior_map` to scale the KL mapped field by a scalar value, 15 in this case. We do this to enable inferring pressure signals with a larger magnitude.

```
G = FEniCSMappedGeometry(G_KL, map=prior_map)
```

We construct a continuous 2D geometry for the observations in which one axis represents the observation times and the other axis represents the sensor locations.

```
G_cont = Continuous2D((obs_times, obs_locations))
```

where `obs_times` and `obs_locations` are arrays of the observation times and locations. Now we create a **CUQIpy** model to encapsulate the forward operator `PAT` with the parametrization, represented by the domain geometry, and the range geometry,

```
A = Model(PAT, domain_geometry=G, range_geometry=G_cont)
```

Note that in creating this model we are treating `PAT` as a black-box function. Now, **CUQIpy-FEniCS** can utilize the information about domain and range geometries to allow advanced sampling and visualization.

The parameterized Bayesian problem for the PAT now takes the form

$$\begin{aligned}\boldsymbol{x} &\sim \mathrm{Gaussian}\left(\boldsymbol{0}, \boldsymbol{I}_{n_{\mathrm{KL}}}\right),\\ \boldsymbol{y} &\sim \mathrm{Gaussian}\left(\boldsymbol{A}(\boldsymbol{x}), s_{\mathrm{noise}}^2 \boldsymbol{I}_{rm}\right).\end{aligned}$$

where $rm$ is the size of $\boldsymbol{y}^{\mathrm{obs}}$. Note that $\boldsymbol{A}$ now contains the KL-expansion parameterization. We can set up this Bayesian problem as

```
x = Gaussian(0, 1, geometry=G)
y = Gaussian(A(x), s_noise**2, geometry=G_cont)
```

In this example we consider the same `s_noise`, as well as, the same noisy data `y_obs` as in section 5.2. Similar to the previous sections we now construct the joint and the posterior distributions.

```
joint = JointDistribution(x, y)
posterior = joint(y=y_obs)
```

Setting up a Bayesian problem for the partial data, e.g. when we place a sensor only on $\xi_L = 0$, is similar. To explore the posterior, we use the preconditioned Crank–Nicolson (pCN) [12] sampler which is suited for the KL parameterization.

```
sampler = pCN(posterior)
samples = sampler.sample_adapt(50000)
```

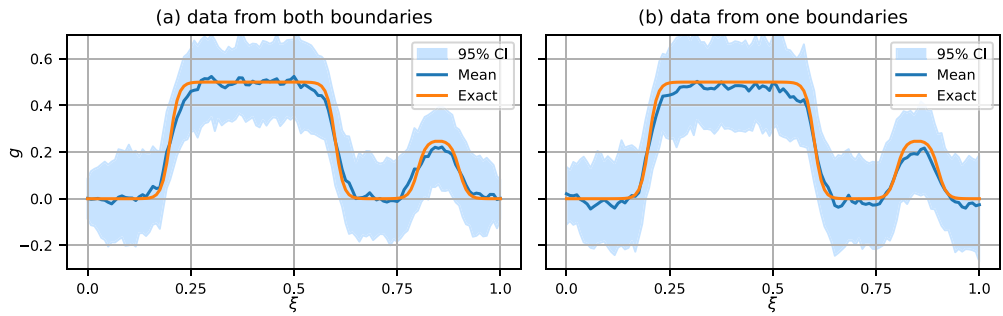

**Figure 11.** Estimated initial pressure profile *g* together with the uncertainty estimates. The plots correspond to the (a) full data and (b) partial data.

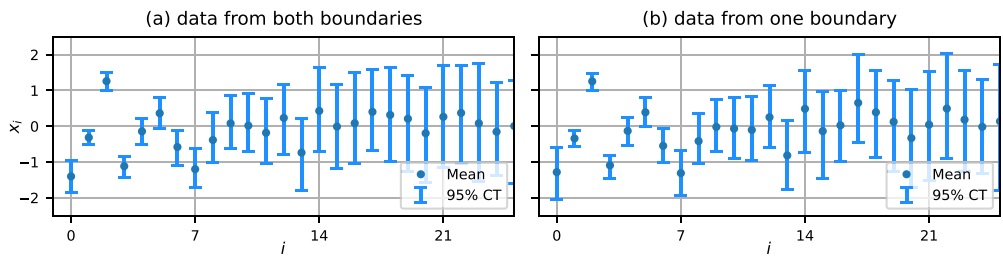

**Figure 12.** Estimation of the first 25 components of ***x***, i.e. KL expansion coefficients in (14), and the uncertainty in this estimation. (a) Full data and (b) partial data.

We can visualize the posterior for ***g***, i.e. the initial pressure distribution, by

```
samples.funvals.vector.plot_ci(95, exact=g_true)
```

This plots the mean function for ***g*** and the CIs associated with this estimate in figure 11. We see that the mean function, in the case with complete data, is a better estimate for the true initial pressure profile compared to the case with partial data. Furthermore, when data corresponding the right boundary is missing, the uncertainty in estimating the right boundary increases.

Finally, we plot the mean and 95% CI for the expansion coefficients ***x***:

```
samples.plot_ci(95, plot_par=True)
```

In figures 12(a) and (b) we present the CI plots for the full data (from both boundaries) and the partial data (only from the left boundary), respectively. Note that we only show the first 25 components although we estimate all 100 parameters. We see that the uncertainty of the first coefficient significantly increases for the case with partial data.

## 6. Conclusion and future work

In this paper we described our general framework for modeling and solving PDE-based Bayesian inverse problems with the CUQIpy Python software package. We showed how to express PDEs natively in CUQIpy, or using a user-provided black-box PDE solver. We also showed how to formulate statistical assumptions about unknown parameters using

CUQIpy and conduct Bayesian inference and uncertainty quantification. We also presented our **CUQIpy-FEniCS** plugin as an example of how to incorporate modeling by third-party PDE libraries such as the finite-element modeling package **FEniCS**.

We showed that **CUQIpy** and **CUQIpy-FEniCS** provide a consistent and intuitive interface to model and solve PDE-based Bayesian inverse problems, as well as analyze and visualize their solutions. Results were shown for parabolic, elliptic and hyperbolic examples involving the heat and Poisson equations as well as application case studies in EIT and PAT.

Future work includes expanding support for derivatives across distributions, forward models, and geometries, as well as integrating **PyTorch** automatic differentiation into **CUQIpy** through the **CUQIpy-PyTorch** plugin. This will simplify the use of gradient-based samplers such as NUTS, as in the Poisson example in this paper, to help address the computational challenge of MCMC-based sampling of high-dimensional and complicated posterior distributions arising in large-scale inverse problems. The extensible plugin structure can also be used to integrate more PDE-based modeling libraries.

Overall, we believe **CUQIpy** and its plugins provide a promising platform for solving PDE-based Bayesian inverse problems and have a significant potential for further development and expansion in the future.

## Data availability statements

**CUQIpy** and plugins are available from https://cuqi-dtu.github.io/CUQIpy. The code and data to reproduce the results and figures of the present paper are available from https://github.com/CUQI-DTU/paper-CUQIpy-2-PDE. The data that support the findings of this study are openly available at the following URL/DOI: https://zenodo.org/doi/10.5281/zenodo.10512535.

## Acknowledgments

This work was supported by The Villum Foundation (Grant No. 25893). J S J would like to thank the Isaac Newton Institute for Mathematical Sciences for support and hospitality during the programme 'Rich and Nonlinear Tomography—a multidisciplinary approach' when work on this paper was undertaken. This work was supported by EPSRC Grant Number EP/R014604/1. This work was partially supported by a grant from the Simons Foundation (J S J). F U has been supported by Academy of Finland (Project Number 353095). The authors are grateful to Kim Knudsen for valuable input in regards to PDE-based inversion in medical imaging, to Aksel Kaastrup Rasmussen for a helpful discussion about the EIT problem formulation and to members of the CUQI project for valuable input that helped shape the design of **CUQIpy**.

## ORCID iDs

Amal M A Alghamdi  https://orcid.org/0000-0003-0145-5296
Nicolai A B Riis  https://orcid.org/0000-0002-6883-9078
Babak M Afkham  https://orcid.org/0000-0003-3203-8874
Felipe Uribe  https://orcid.org/0000-0002-1010-8184
Silja L Christensen  https://orcid.org/0000-0003-3995-3055
Per Christian Hansen  https://orcid.org/0000-0002-7333-7216
Jakob S Jørgensen  https://orcid.org/0000-0001-9114-754X